\documentclass{article}
\usepackage{arxiv}
\usepackage{indentfirst}
\usepackage{amsmath}
\usepackage[utf8]{inputenc} 
\usepackage[T1]{fontenc}    
\usepackage{hyperref}       
\usepackage{url}            
\usepackage{booktabs}       
\usepackage{amsfonts}       
\usepackage{nicefrac}       
\usepackage{microtype}      
\usepackage{graphicx}
\usepackage{float}
\usepackage{afterpage}
\usepackage{booktabs, calc}
\usepackage{multicol}
\usepackage{multirow}
\usepackage{indentfirst}
\usepackage{algorithm}
\usepackage{algpseudocode}
\usepackage{subcaption}
\usepackage{algorithm} 

\DeclareMathAlphabet\mathbfcal{OMS}{cmsy}{b}{n} 
\usepackage{lineno}

\newtheorem{theorem}{Theorem}

\newtheorem{definition}[theorem]{Definition}

\newcommand{\R}{\mathbb R}

\renewcommand{\Re}{\mathrm{Re} \;}

\newcommand\norm[1]{\left\lVert#1\right\rVert}

\title{Efficiently and accurately simulating multi-dimensional M-coupled nonlinear Schrödinger equations with fourth-order time integrators and Fourier spectral method}

\author{
Nathaniel Lovett  \\
Department of Mathematics\\
Utah Valley University\\
Orem, UT 84058\\
  \texttt{nlovett917@gmail.com} \\
   \And
Harish Bhatt \\
  Department of Mathematics\\
  Utah Valley University\\
 Orem, UT 84058 \\
  \texttt{harish.bhatt@uvu.edu} \\
}

\begin{document}
\maketitle
\begin{abstract}
Coupled nonlinear Schrödinger equations model various physical phenomena, such as wave propagation in nonlinear optics, multi-component Bose-Einstein condensates, and shallow water waves. Despite their extensive applications, analytical solutions of coupled nonlinear Schrödinger equations are widely either unknown or challenging to compute, prompting the need for stable and efficient numerical methods to understand the nonlinear phenomenon and complex dynamics of systems governed by coupled nonlinear Schrödinger equations. This paper explores the use of the fourth-order Runge-Kutta based exponential time-differencing and integrating factor methods combined with the Fourier spectral method to simulate multi-dimensional M-coupled nonlinear Schrödinger equations.

The theoretical derivation and stability of the methods, as well as the runtime complexity of the algorithms used for their implementation, are examined. Numerical experiments are performed on systems of two and four multi-dimensional coupled nonlinear Schrödinger equations. It is demonstrated by the results that both methods effectively conserve mass and energy while maintaining fourth-order temporal and spectral spatial convergence. Overall, it is shown by the numerical results that the exponential time-differencing method outperforms the integrating factor method in this application, and both may be considered further in modeling more nonlinear dynamics in future work.
\end{abstract}


\section{Introduction}
\label{sec:1}
\indent The system of coupled nonlinear Schrödinger equations (CNLSE), introduced in 1967, models the propagation and interaction of multiple wave components in a nonlinear medium \cite{BN}. These equations are applicable in various fields such as nonlinear optics and transmission lines \cite{Yang, WLEK}, water wave dynamics \cite{DLS, WLEK}, Bose-Einstein condensates \cite{ID, WLEK}, plasma physics \cite{Sulem}, and quantum mechanics \cite{FLMN, CC}. In this work, we consider the following multi-dimensional M-coupled nonlinear Schrödinger equation system with no external potential:
\begin{flalign}
&i\frac{\partial \Psi_j(x, t)}{\partial t} - \alpha_j(-\Delta)\Psi_j(x,t) + \left(\sum_{m=1}^{M}\sigma_{jm}|\Psi_m(x,t)|^2\right)\Psi_j(x,t) = 0, & \label{eq:1} \\
&x\in\Omega\subset\mathbb{R}^d, \quad t>0, \quad j = 1,2,...,M & \notag
\end{flalign}
with periodic, homogeneous Neumann, or Dirichlet boundary conditions, and initial conditions:
\begin{flalign*}
&\Psi_j(x,0)=\Psi_{j_0}(x), \quad x\in\Omega,&
\end{flalign*}
where $\Psi_j$ are complex wave amplitudes, $i=\sqrt{-1}$, $x$ and $t$ are space and time variables respectively, $\sigma_{jm}$ are self-phase modulation coefficients for $j=m$ and cross-phase modulation coefficients for $j \neq m$, $\alpha_j$ are the group velocity dispersion coefficients, $-\Delta$ represents the multi-dimensional Laplacian, $M$ represents the number of coupled equations, and $\Omega \subset \mathbb{R}^d, d=1,2,3$ is the bounded domain with Lipschitz continuous boundary in $\mathbb{R}^d$.

There are a variety of specific cases where CNLSE systems have analytical solutions. When $M=1$, \eqref{eq:1} becomes the well-known integrable cubic NLSE \cite{WLEK}. For $M=2$, with $\alpha_n=1$ and $\sigma_{jj}=e$, the CNLSE system becomes the integrable Manakov equations, which can be solved exactly using the inverse scattering method \cite{BK}. For example, in \cite{WIH}, Wadati et al. exactly solved CNLSEs for $M=2$ that models linear birefringence in nonlinear optical fibers under Dirichlet boundary conditions. Yang investigated this class of analytical solutions further in \cite{Yang}. Lin and Wei proved conditions for the existence and nonexistence of ground-state solutions of steady-state CNLSE systems in \cite{LW}. Moreover, Seadawy and Cheema found a variety of solitary wave solutions to CNLSEs with Kerr law nonlinearity \cite{SC}. However, many CNLSE systems lack analytical solutions, especially when they are non-integrable, multi-dimensional, or model chaotic systems with complex initial conditions \cite{WLEK}.

To address CNLSE systems without known analytical solutions, extensive research has focused on developing efficient and stable numerical methods to assist in approximating possible CNLSE solutions, such as finite difference, finite element, operator splitting, exponential time differencing (ETD), integrating factor (IF), spectral, and pseudo-spectral methods. In \cite{MIT}, Ismail and Taha presented a stable second-order finite difference method for CNLSE systems modeling solitons in linearly birefringent optical fibers. In \cite{ISM}, a fourth-order finite difference method is utilized with explicit RK4 methods, achieving fast computation with conditional stability. Alamri et al. developed a highly accurate finite difference method using Newton's method and a predictor-corrector method, finding Newton's method more effective under their method in \cite{IA}. In \cite{Ashi}, Ashi compared ETD and IF methods in combination with Fourier spectral approximation for modeling several different nonlinear and semilinear PDEs, including NLSEs. They concluded that ETD methods are better for non-traveling or slowly traveling soliton solutions, while IF methods are preferable for solitons solutions with a large speed. A local discontinuous Galerkin (LDG) finite element method in combination with a fourth-order Runge-Kutta ETD method (ETDRK4) was introduced by Liang et al. in \cite{LKX} to solve n-coupled 1D CNLSE systems that is shown to be efficient, stable, and conserves energy and mass over a long time scale. Qian et al. in \cite{QSC}, introduced a semi-explicit multi-symplectic splitting method to solve a 3-coupled NLSE system, finding it to be more efficient than a known reliable implicit multi-symplectic splitting method. In \cite{HBK}, Bhatt and Khaliq implement an improved version of Cox and Matthews ETD3RK method in combination with the Numerov/Douglas approximation to solve \eqref{eq:1} for $M=2,4$, concluding the method to be efficient and stable in approximating long range solitary solutions. Recently in \cite{musta}, Almushaira developed several conservative higher-order compact finite difference methods to solve the two-dimensional (2D) NLSE, finding them to be accurate and stable with expected orders of convergence.

In this manuscript, we introduce two fourth-order Runge-Kutta-based time-stepping methods, combined with the Fourier spectral method, to approximate the spatio-temporal evolution of system \eqref{eq:1}. One method is a modification of Krogstad's ETDRK4-B method, and the other is a modification of the Integrating Factor Runge-Kutta method (IFRK4). These numerical methods have been successfully applied to solving various nonlinear PDEs, as noted in \cite{Cox, HBK, JLQY}. However, to the best of our knowledge, these methods have not been used to approximate solutions for the 4-coupled NLSE system in higher dimensions. Therefore, our primary focus in this work is to analyze the effectiveness of these methods in combination with the Fourier spectral method for approximating solutions to system \eqref{eq:1} and to evaluate their conservation and stability properties.

The remainder of this paper is organized as follows. In Section \ref{sec:2}, we discuss the Fourier spectral method used to discretize the spatial derivatives and transform the coupled nonlinear Schrödinger equations (CNLSEs) into ordinary differential equations (ODEs) in Fourier frequency space. Section \ref{sec:3} introduces two time-stepping methods for solving the resulting ODEs, along with an analysis of their linear stability. In Section \ref{sec:4}, we examine conserved wave properties to evaluate the applicability of the proposed methods. We then conduct several numerical experiments on multi-dimensional CNLSEs under various boundary conditions to assess the performance of the proposed methods in terms of accuracy and computational efficiency. In Section \ref{sec:5}, we briefly conclude our findings.

\section{Spatial discretization}
\label{sec:2}
\indent The Fourier spectral method is an effective technique for solving multidimensional PDEs. It leverages the Fourier transform's properties to handle spatial derivatives efficiently. In this section, we describe how the Fourier spectral method can be applied to solve multidimensional CNLSEs involving the Laplacian operator $-\Delta$. Following the derivation in \cite{BOKB}, we suppose there is a complete set of orthonormal trigonometric eigenfunctions $\phi_{j_{\gamma_1, \ldots, \gamma_d}}$ with corresponding eigenvalues $\lambda_{j_{\gamma_1, \ldots, \gamma_d}}$ for $j=1,\ldots,M$ satisfying periodic, homogeneous Neumann, or homogeneous Dirichlet boundary conditions for the Laplacian operator on the bounded $d$-dimensional domain $\Omega = [a, b]^d \subset \mathbb{R}^d$:
\begin{flalign*}
& (-\Delta) \phi_{j_{\gamma_1,\ldots,\gamma_d}} = \lambda_{j_{\gamma_1,\ldots,\gamma_d}} \phi_{j_{\gamma_1,\ldots,\gamma_d}}.&
\end{flalign*}
Given this, the solution \(\Psi_j\) can be expressed as:
\begin{flalign*}
&\Psi_{j_{\gamma_1, \ldots, \gamma_d}} = \sum_{\gamma_1=0}^{\infty} \cdots \sum_{\gamma_d=0}^{\infty} \hat{\Psi}_{j_{ \gamma_1,\ldots, \gamma_d}} \phi_{j_{\gamma_1,\ldots,\gamma_d}}&
\label{eq:Fourierinf}
\end{flalign*}
where the spectral coefficients \(\hat{\Psi}_{j_{\gamma_1,\ldots, \gamma_d}}\) are defined as the inner product of \(\Psi_j\) with the eigenfunctions \(\phi_{j_{\gamma_1,\ldots,\gamma_d}}\):
\begin{flalign*}
&\hat{\Psi}_{j_{\gamma_1,\ldots, \gamma_d}} = \langle \Psi_j, \phi_{j_{\gamma_1,\ldots,\gamma_d}} \rangle = \int_{\Omega} \Psi_j(x) \phi_{j_{\gamma_1,\ldots,\gamma_d}}(x) \, dx &
\end{flalign*}
Efficient computation of the coefficients \(\hat{\Psi}_{j_{\gamma_1,\ldots, \gamma_d}}\) and the inverse reconstruction of \(\Psi_j\) in physical space can be achieved using fast and robust algorithms such as the direct and inverse Discrete Sine Transforms (DST) for homogeneous Dirichlet boundary, the direct and inverse Discrete Cosine Transforms (DCT) for homogeneous Neumann boundary, and the direct and inverse Fast Fourier Transforms (FFT) for periodic boundary conditions \cite{BOKB, BH, BOPG}. In the context of the Fourier spectral method, \(-\Delta \Psi_j\) can be approximated using N eigenfunctions, where N is the number of internal equispaced spatial discretization points:
\begin{flalign}
&-\Delta \Psi_j \approx \sum_{\gamma_1=0}^{N-1} \cdots \sum_{\gamma_d=0}^{N-1} \hat{\Psi}_{j_{\gamma_1,\ldots, \gamma_d}} \lambda_{j_{\gamma_1,\ldots,\gamma_d} }\phi_{j_{\gamma_1,\ldots,\gamma_d}}&
\end{flalign}
This spectral decomposition is well known and has been used to find analytical solutions to PDE problems \cite{BOKB}. The different boundary conditions we consider determine the specific spectral decomposition of the Laplacian operator, with each utilizing different eigenvalues, eigenfunctions, and equispaced mesh grids \cite{BOKB, BH, BOPG}. With $L=b-a$ and $X=\{x_i\} \in \Omega^d$, the eigenfunctions, eigenvalues, and mesh grids required under several boundary conditions are defined as below:
\paragraph{Periodic Boundary Conditions}
For periodic boundaries, the eigenvalues, eigenfunctions, and mesh grid are given by:
\begin{flalign*}
\lambda_{j_{\gamma_1, \ldots, \gamma_d}} &= \sum_{i=1}^{d} \left( \frac{2\pi \gamma_i}{L} \right)^2 \ \forall \ \gamma_1\ldots\gamma_d, &\\
\phi_{j_{\gamma_1, \ldots, \gamma_d}} &= \prod_{i=1}^{d} \left(\sqrt{\frac{1}{L}}\right) e^{i 2\pi \gamma_i (X-a)/L} \ \forall \ \gamma_1\ldots\gamma_d,&
\\ x_i &= -\frac{L}{2} + (i)h, \ i=1,\ldots, N \ , \quad h=\frac{L}{N}.&
\end{flalign*}

\paragraph{Homogeneous Dirichlet Boundary Conditions}
Under these conditions, the parameters become:
\begin{flalign*}
\lambda_{j_{\gamma_1, \ldots, \gamma_d}} &= \sum_{i=1}^{d} \left( \frac{(\gamma_i + 1)\pi}{L} \right)^2 \ \forall \ \gamma_1\ldots\gamma_d &\\
\phi_{j_{\gamma_1, \ldots, \gamma_d}}, &= \prod_{i=1}^{d} \left(\sqrt{\frac{2}{L}}\right) \sin\left( \frac{(\gamma_i + 1)\pi (X - a)}{L} \right) \ \forall \ \gamma_1\ldots\gamma_d,&
\\ x_i &= a + (i)h,  \ i=1,\ldots, N \ , \quad h=\frac{L}{N+1}. &
\end{flalign*}

\paragraph{Homogeneous Neumann Boundary Conditions}
For these conditions, the parameters become:
\begin{flalign*}
\lambda_{j_{\gamma_1, \ldots, \gamma_d}} &= \sum_{i=1}^{d} \left( \frac{\gamma_i \pi}{L} \right)^2 \ \forall \ \gamma_1\ldots\gamma_d, &\\
\phi_{j_{\gamma_1, \ldots, \gamma_d}} &= \prod_{i=1}^{d} \sqrt{\frac{2}{L}} \cos\left( \frac{\gamma_i (X-a) \pi}{L} \right) \ \forall \ \gamma_1\ldots\gamma_d,&
\\ x_i &= a + (i-1)h + \frac{h}{2},  \ i=1,\ldots, N \ , \quad h=\frac{L}{N}.&
\end{flalign*}

To demonstrate how the Fourier spectral method applied to the Laplacian, let's consider its application to the 1D CNLSE for simplicity. The 1D CNLSE in Fourier space using the Fourier spectral technique and applying the Fourier transform $\mathcal{F}$ is represented as follows:

\begin{flalign} 
&\frac{d \mathbf{\hat{\Psi}}_j(t)}{dt} + A(\mathbf{\hat{\Psi}}_j) = \mathbf{\hat{F}}(\Psi_j), &
\label{eq:transformedeq} \\
&t \in (0,T], \quad \mathbf{\hat{\Psi}}_j(0)= \mathbf{\hat{\Psi}}_{j_0}, \quad j=1,\ldots,M \notag &
\end{flalign}
where $A=i\alpha_j(\lambda_{j_{\gamma_{1_1}}},\ldots,\lambda_{j_{\gamma_{1_N}}})$, $\mathbf{\hat{\Psi}}_j := \mathcal{F}([\Psi_{j}(x_1,t), \Psi_{j}(x_2,t), \ldots, \Psi_{j}(x_N,t)]^{T})$, and $\mathbf{\hat{F}}(\Psi_j) = i\mathcal{F}([f_j(x_1, t, \Psi_j), f_j(x_2, t, \Psi_j), \ldots, f_j(x_N, t, \Psi_j)]^{T})$ where $f_j(x_i,t,\Psi_j)$ is a function of the nonlinear terms in the $j$th Equation of system \eqref{eq:1}. In practical computations, the mesh choice will depend on the specified boundary conditions, and grids can be generated in MATLAB as follows (see \cite{Tref} for more information): 

\noindent For periodic boundary conditions,
\begin{flalign*}
& x = \frac{L}{N}\left[-\frac{N}{2}+1:\frac{N}{2}\right], &\\
& \lambda_{j_{\gamma_1}} = \left(\frac{2\pi}{L}\big[0:\frac{N}{2} ~~ -\frac{N}{2}+1: -1\big]\right)^2. &
\end{flalign*}

\noindent For homogeneous Dirichlet boundary conditions,
\begin{flalign*}
& x = a + [1:N]h, \quad h = \frac{L}{N+1}, &\\
& \lambda_{j_{\gamma_1}} = \left(\frac{\pi}{L}[1:N]\right)^2. &
\end{flalign*}

\noindent For homogeneous Neumann boundary conditions,
\begin{flalign*}
& x = a + [0:N-1]h + \frac{h}{2}, \quad h = \frac{L}{N}, &\\
& \lambda_{j_{\gamma_1}} = \left(\frac{\pi}{L}[0:N-1]\right)^2. &
\end{flalign*}

While we consider a one-dimensional (1D) space for simplicity, the higher dimensional cases are defined analogously.

\section{Time stepping methods and solution strategy}
\label{sec:3}
\subsection{Fourth-Order Exponential Time Differencing Runge-Kutta Method}
\indent In this subsection, we discuss the formation of ETD methods and the Krogstad-P22 method we consider for this work. ETD methods have received more attention and treatment in recent years in approximating nonlinear PDEs \cite{HIRS, BK, HBK, KRO, LKX}, particularly those that suffer from stiffness \cite{Cox}. To briefly describe the ETD method, here we consider system \eqref{eq:1} in 1D and note the parameters are defined analogously in higher dimensions. Because we utilize the Fast Fourier Transform (FFT) to implement the Fourier spectral method computationally, we use the proper inverse fast Fourier transform (IFFT) to reacquire the solution in physical space.

By discretizing $-\alpha_j(-\Delta)\Psi_j$ in \eqref{eq:1} with the Fourier spectral method as shown in Section \ref{sec:2}, we acquire a system of coupled ODEs in Fourier space shown in Equation \eqref{eq:transformedeq} with appropriate boundary conditions. Let $k=t_{n+1}-t_n$ be the time step size then we come with the following formula for $\mathbf{\hat{\Psi}}_j(t_{n+1})=\mathbf{\hat{\Psi}}_{j_{n+1}}$ utilizing the variation of constants formula:
\begin{flalign}
    &\hat{\mathbf{\Psi}}_{j_{n+1}}=e^{-kA}\hat{\mathbf{\Psi}}_{j_n}+k \int_{0}^{1}e^{-kA(1-\tau)}\hat{\mathbf{F}}(\mathbf{\Psi}_j(x,t_n+\tau k),t_n+\tau k)d\tau.
\label{eq:etdint}&
\end{flalign}
\indent There are many ways to approximate the integral and matrix exponential terms in Equation \eqref{eq:etdint}. These choices of approximations determine the specific ETD method in use. Cox and Matthews developed the ETDRK4 method in \cite{Cox} that utilizes a fourth-order Runge-Kutta approximation of the integral in \eqref{eq:etdint}. This method has been of particular interest due to the accuracy and efficiency in handling nonlinear differential equations \cite{LK, KRO, LKX, CK, BK}. However, ETDRK4 is known to suffer from numerical instability in situations of small or zero eigenvalues of the linear operator $A$ \cite{Cox}. Kassam and Trefethen sought to handle this by solving the integral in \eqref{eq:etdint} as a contour integral over a path that contains the spectrum of $A$ \cite{KT}. As discussed in \cite{HBK}, this approach becomes difficult when $A$ has larger spectral radii and is computationally inefficient for larger problems in multiple dimensions. Krogstad refined the ETDRK4 method alternatively in \cite{KRO}, improving the accuracy and stability by utilizing different choices of $a_n$, $b_n$, and $c_n$ derived from commutator-free Lie group methods. Krogstad's refined method is denoted as the ETDRK4-B method. The ETD method we consider in this work is a modified ETDRK4-B method, that uses a Pad$\acute{\mathrm{e}}$-(2,2) approximation for matrix exponential terms, that was introduced by Bhatt and Khaliq for approximating solutions to coupled Burgers' equations in \cite{HBK}. This Pad$\acute{\mathrm{e}}$-(2,2) approximation has the advantage of cancellation of $A^{-1}$ and $A^{-3}$ terms, reducing computations needed.
\\ \indent We use the notation $R_{r,s}(z)$ for $(r, s)-$Pad$\acute{\mathrm{e}}$ approximation to $e^{-z}$. The $(r+s)^{th}$ order rational Pad$\acute{\mathrm{e}}$ approximation to $e^{-z}$ is defined as:
\begin{flalign*}
&R_{r,s}(z)=\frac{P_{r,s}(z)}{Q_{r,s}(z)}, &
\end{flalign*}
where
\begin{flalign*}
&P_{r,s}(z)=\sum\limits_{j=0}^r\frac{(s+r-j)!r!}{(s+r)!j!(r-j)!}(-z)^j\ \text{and}\ 
Q_{r,s}(z)=\sum\limits_{j=0}^s\frac{(s+r-j)!}{(s+r)!j!(s-j)!}(z)^j. &\\
\end{flalign*}

\begin{definition}[Pad$\acute{\mathrm{e}}$-(r,s) Rational Approximation] A rational approximation $R_{r,s}(z)$ of $e^{-z}$ is said to be $A-$acceptable, if $|R_{r,s}(z)|<1$, whenever $\Re(z)<0$ and $L-$acceptable if, in addition, $|R_{r,s}(z)|\rightarrow0\ \text{as}\ \Re(z)\rightarrow-\infty.$\\
\indent The rational Pad$\acute{\mathrm{e}}$ approximation $R_{r,s}(z)$ to $e^{-z}$ is :
\begin{itemize}
\item $A-$acceptable if $r=s$.
\item $L-$acceptable if $r=s-1$ or $s-2$.
\end{itemize}
\end{definition}
When we utilize the A-acceptable (2,2)-Pad$\acute{\mathrm{e}}$ approximation of $e^{-z}$ given by $R_{2,2}(z) = \frac{12-6z+z^{2}}{12+6z+z^{2}}$ in the ETDRK4-B method from \cite{KRO}, we get the following formula as given in \cite{HBK}:
\begin{flalign}
    &\mathbf{\hat{\Psi}}_{j_{n+1}}=R_{2,2}(kA)\mathbf{\hat{\Psi}}_{j_n}+P_1(kA)F_n+P_2(kA)(-3F_n+2F_{n}^{a}+2F_{n}^{b}-F_{n}^{c})+P_3(kA)(F_n-F_{n}^{a}-F_{n}^{b}+F_{n}^{c}), &
\label{eq:krogp22}
\end{flalign}
where
\begin{flalign}
    &F_n=\mathbf{\hat{F}}(\mathbf{\Psi}_{j_n},t_n), \quad F_{n}^{a}=\mathbf{\hat{F}}(\textstyle a_n,t_n+\frac{k}{2}), &\notag \\
    &F_{n}^{b}=\mathbf{\hat{F}}(\textstyle b_n,t_n+\frac{k}{2}), \quad F_{n}^{c}=\mathbf{\hat{F}}(c_n,t_n+k), &\notag \\
    &P_1(kA)=12k(12+6kA+k^2A^2)^{-1}, &\notag \\
    &P_2(kA)=k(6-kA)(12+6kA+k^2A^2)^{-1}, &\notag \\
    &P_3(kA)=2k(4-kA)(12+6kA+k^2A^2)^{-1}, &\label{eq:FnPn}
\end{flalign}
and
\begin{flalign}
    &a_n=\mathcal{F}^{-1}(\Tilde{R}_{2,2}(kA)\mathbf{\hat{\Psi}}_{j_n}+\Tilde{P_1}(kA)F_n), &\notag \\
    &b_n=\mathcal{F}^{-1}(\Tilde{R}_{2,2}(kA)\mathbf{\hat{\Psi}}_{j_n}+\Tilde{P_1}(kA)F_n+\Tilde{P_2}(kA)(F_{n}^{a}-F_n)), &\notag \\
    &c_n=\mathcal{F}^{-1}(R_{2,2}(kA)\mathbf{\hat{\Psi}}_{j_n}+P_1(kA)F_n+2P_2(kA)(F_{n}^{b}-F_n)), &\label{eq:rkcoeffs}
\end{flalign}
such that
\begin{flalign}
    &R_{2,2}(kA)=(12-6kA+k^2A^2)(12+6kA+k^2A^2)^{-1}, &\notag \\
    &\Tilde{R}_{2,2}(kA)=(48-12kA+k^2A^2)(48+12kA+k^2A^2)^{-1}, &\notag \\
    &\Tilde{P_1}(kA)=24k(48+12kA+k^2A^2)^{-1}, &\notag \\
    &\Tilde{P_2}(kA)=2k(12+kA)(48+12kA+k^2A^2)^{-1}. &\label{hatp22}
\end{flalign}

We will refer to this method as Krogstad-P22. Below, in Algorithm \ref{alg:krogalg}, we show the FFT implementation of the Krogstad-P22 method in three dimensions (3D) with homogeneous Neumann boundary conditions.
\begin{algorithm}[H]
\caption{3D Krogstad-P22 with homogeneous Neumann boundary conditions}
\label{alg:krogalg}
\begin{algorithmic}[1]
    \State Given $L= x_n-x_0,\ N, h = \frac{L}{N}, \alpha, T$, and $k$
    \State Compute $\lambda = \Big((0:N-1)\frac{\pi}{L}\Big)$;~$\lambda_y = \lambda_z = \lambda_x$; ~ $x = x_0+(0:N-1)h+\frac{h}{2}$;~~$y=x$;~~$z=x$;
    \Statex $ A = i \alpha (\lambda_x^2+\lambda_y^2+\lambda_z^2) $;~~ $\Psi_j(x,y,z,0)=\Psi_{j_0}$;~$t=[0:k:T]$; ~~$m=\text{len}(t)$
    \Statex  Precompute the following arrays:
    \State $P_1(kA)=12k(12+6kA+k^2A^2)^{-1}$
    \State $P_2(kA)=k(6-kA)(12+6kA+k^2A^2)^{-1}$
    \State $P_3(kA)=2k(4-kA)(12+6kA+k^2A^2)^{-1}$
    \State $R_{2,2}(kA)=(12-6kA+k^2A^2)(12+6kA+k^2A^2)^{-1}$
    \State $\Tilde{R}_{2,2}(kA)=(48-12kA+k^2A^2)(48+12kA+k^2A^2)^{-1}$
    \State $\Tilde{P_1}(kA)=24k(48+12kA+k^2A^2)^{-1}$
    \State $\Tilde{P_2}(kA)=2k(12+kA)(48+12kA+k^2A^2)^{-1}$
    \For{$ n=0,\cdots, m-1$} 
    \State $\hat{\Psi}_j,n=dctn(\Psi_j,n)$
    \Comment{\textit{dctn} is MATLAB built-in function} 
    \State $a_n=\textit{idctn}(\Tilde{R}_{2,2}(kA)\hat{\Psi}_{j_n}+\Tilde{P}_1(kA)F_n$
    \Comment{\textit{idctn} is MATLAB built-in function}
    \State $b_n=\textit{idctn}(\Tilde{R}_{2,2}(kA)\hat{\Psi}_{j_n}+ \Tilde{P}_1(kA)F_n+ \Tilde{P}_2(kA)(F_n^a-F_n))$
    \State $c_n=\textit{idctn}(R_{2,2}(kA)\hat{\Psi}_{j_n}+ P_1(kA)F_n + 2 P_2(kA)(F_n^b-F_n))$
    \State $\Psi_{j_{n+1}} = \textit{idctn}(R_{2,2}(kA)\hat{\Psi}_{j_n} + F_n(P_1(kA)-3 P_2(kA)+P_3(kA)) + F_n^a(2 P_2(kA)-P_3(kA))$
    \Statex \quad \quad \quad \quad \quad \quad \quad \quad$+ F_n^b(2 P_2(kA)-P_3(kA)) + F_n^c(P_3(kA)-P_2(kA)))$
    \EndFor
\end{algorithmic}
\end{algorithm}

\subsection{Fourth-Order Integrating Factor Runge-Kutta method}
\indent In this subsection, we briefly describe the formation of IF methods and the IFRK4 method given in \cite{JLQY}, and modify this method to make it more efficient for solving the multi-dimensional system \eqref{eq:1}. Integrating factor methods utilize the change of variables $\mathbf{w}_j=\mathbf{\hat{\Psi}}_je^{kA}$ after multiplying system \eqref{eq:transformedeq} by the integrating factor $e^{kA}$, then approximating the changed system with a numerical method before changing variables back to $\mathbf{\hat{\Psi}}_j$. This choice of numerical method determines the specific IF method in use \cite{Cox,KT,JLQY}.
In \cite{JLQY}, Ju et al. use the fourth-order Runge-Kutta method to solve the system in $\mathbf{w}_j$, resulting in the following IFRK4 method when variables are changed back to $\mathbf{\hat{\Psi}}_j$:
\begin{flalign}
    &\mathbf{\hat{\Psi}}_{j_{n+1}}= e^{-kA}\mathbf{\hat{\Psi}}_{j_{n+1}}+ k\Big({\textstyle \frac{1}{6}}e^{-kA}F_n+{\textstyle \frac{1}{3}}e^{\frac{-kA}{2}}F^a_n+{\textstyle \frac{1}{3}}e^{\frac{-kA}{2}}F^b_n+{\textstyle \frac{1}{6}}F^c_n\Big), & \\
    &\text{where} \notag& \\
    &a_n=\mathcal{F}^{-1}\Big(e^{\frac{-kA}{2}}(\mathbf{\hat{\Psi}}_{j_n}+{\textstyle \frac{-k}{2}}F_n)\Big), \notag& \\
    &b_n=\mathcal{F}^{-1}\Big(e^{\frac{-kA}{2}}\mathbf{\hat{\Psi}}_{j_n}+{\textstyle \frac{-k}{2}F^a_n}\Big), \notag& \\
    &c_n=\mathcal{F}^{-1}\Big(e^{-kA}\mathbf{\hat{\Psi}}_{j_n}+ke^{\frac{-kA}{2}}F^c_n\Big), \notag& \\
    &\text{and} \notag& \\
    &F_n = \mathbf{\hat{F}}(\mathbf{\Psi}_{j_n},t_n), F_{n}^{a}=\mathbf{\hat{F}}(a_n,t_n+ \textstyle \frac{k}{2}), &\notag \\
    &F_{n}^{b}=\mathbf{\hat{F}}(b_n,t_n+ \textstyle \frac{k}{2}), F_{n}^{c}=\mathbf{\hat{F}}(c_n,t_n+k). &\notag
\end{flalign}

In MATLAB, the built-in \textit{expm} function, which has order of complexity $\mathcal{O}(N^3)$, is used to handle the matrix exponential terms efficiently in 1D. However, for higher dimensions or a larger system in \eqref{eq:1}, using \textit{expm} becomes computationally costly. There are many methods that one can consider for handling matrix exponential terms \cite{ML}. For efficient computation in higher dimensions, we modify this method by approximating the matrix exponential term $e^{kA}$ with an L-acceptable (1,3)-Pad$\acute{\mathrm{e}}$ approximation of $e^{-z}$ given by $R_{1,3}(z) = \frac{(1-\frac{z}{4})}{(1+\frac{3}{4}z+\frac{z^2}{4}+\frac{z^3}{24})}$, which yields the following modified solution for $\hat{\Psi}_{j_{n+1}}$:
\begin{flalign}
    &\mathbf{\hat{\Psi}}_{j_{n+1}}=R_{1,3}(kA) \mathbf{\hat{\Psi}}_{j_n}+k \Big( {\textstyle \frac{1}{6}}R_{1,3}(kA)F_n+{\textstyle \frac{1}{3}}\Tilde{R}_{1,3}(kA)F_{n}^{a}+{\textstyle \frac{1}{3}}\Tilde{R}_{1,3}(kA)F_{n}^{b}+{\textstyle \frac{1}{6}}F_{n}^{c} \Big),&
\label{eq:ifrk}
\end{flalign}
where
\begin{flalign}
    &a_n=\mathcal{F}^{-1}\left(\Tilde{R}_{1,3}(kA) \left(\mathbf{\hat{\Psi}}_{j_n}+\frac{k}{2}F_n\right)\right), &\notag \\
    &b_n=\mathcal{F}^{-1}\left(\Tilde{R}_{1,3}(kA) \mathbf{\hat{\Psi}}_{j_n}+\frac{k}{2}F_{n}^{a}\right), &\notag \\
    &c_n=\mathcal{F}^{-1}\left(R_{1,3}(kA)\mathbf{\hat{\Psi}}_{j_n}+k\Tilde{R}_{1,3}(kA)F_{n}^{b}\right), &\notag 
\end{flalign}
such that
\begin{flalign}
    &R_{1,3}(kA)=(24-6kA)(24+18kA+6k^2A^2+k^3A^3)^{-1}, &\notag \\
    &\Tilde{R}_{1,3}(kA)=24(8-kA)(192+72kA+12k^2A^2+k^3A^3)^{-1}. &\notag
\end{flalign}

We will refer to this modified method as IFRK4-P13. For the Pad$\acute{\mathrm{e}}$ approximations used in Krogstad-P22 and IFRK4-P13, there are limitations in their feasible application. The matrices needed are computed and stored during initialization of computation to save computing time, but calculating inverses of higher order matrix polynomials can yield instability and high round-off error for matrices A with high condition numbers. This can be alleviated by utilizing a partial fraction decomposition as noted in \cite{LK,BK,HBK,KLF}. However, an advantage in using the Fourier Spectral Method for discretization of the spatial domain is that the linear operator $A$ in \eqref{eq:transformedeq} is diagonal in the frequency domain, thus alleviating the need for further partial fraction decomposition. Below in Algorithm \ref{alg:ifrkalg}, we show the implementation of the IFRK4-P13 method in 3D with homogeneous Dirichlet boundary conditions.
We can see how the defined parameters of \eqref{eq:transformedeq} are analogously defined in the higher dimension.
\begin{algorithm}[H]
\caption{3D IFRK4-P13 with homogeneous Dirichlet boundary conditions}
\label{alg:ifrkalg}
\begin{algorithmic}[1]
    \State Given $L= x_n-x_0, N, h = \frac{L}{N+1}, \alpha, T$, $k$.
    \State Compute $\lambda_{x} = \Big((1:N)\frac{\pi}{L}\Big)$;~~$\lambda_y = \lambda_z = \lambda_x$; $x = x_0+(1:N)h$~;~$y=x$;~$z=x$;
    \Statex $ A = i \alpha (\lambda_x^2+\lambda_y^2+\lambda_z^2) $;~~ $\Psi_j(x,y,z,0)=\Psi_{j_0}$;~$t=[0:k:T]$; ~~$m=len(t)$
    \Statex Precompute the following arrays:
    \State $R_{1,3}(kA)=(24-6kA)(24+18kA+6k^2A^2+k^3A^3)^{-1}$
    \State $\Tilde{R}_{1,3}(kA)=24(8-kA)(192+72kA+12k^2A^2+k^3A^3)^{-1}$
    \For{$ n=0,\cdots, m-1$} 
    \State $\hat{\Psi}_{j_n}=dstn(\Psi_{j_n})$
    \Comment{\textit{dstn} is MATLAB built-in function} 
    \State $a_n = \textit{idstn}(\Tilde{R}_{1,3}(kA)(\hat{\Psi}_{j_n}+\frac{k}{2}F_n))$
    \Comment{\textit{idstn} is MATLAB built-in function}
    \State $b_n = \textit{idstn}(\Tilde{R}_{1,3}(kA)\hat{\Psi}_{j_n} + \frac{k}{2}F_n^a)$ 
    \State $c_n = \textit{idstn}(R_{1,3}(kA)\hat{\Psi}_{j_n} + k \Tilde{R}_{1,3}(kA) F_n^b)$
    \State $\Psi_{j_{n+1}} = \textit{idstn}(R_{1,3}(kA)\hat{\Psi}_{j_n} + k (\frac{1}{6} R_{1,3}(kA) F_n + \frac{1}{3} \Tilde{R}_{1,3}(kA) (F_n^a+F_n^b) + \frac{1}{6}F_n^c))$ 
    \EndFor
\end{algorithmic}
\end{algorithm}

\subsection{Time Complexity Analysis}
\indent We briefly consider the time complexity of Krogstad-P22 and IFRK4-P13 algorithms as we have described them. To understand the runtime complexity of these algorithms, we conduct basic Big-O analysis on them with respect to a single iteration of each method, utilizing the following definition:

\begin{definition}[Big-$\mathcal{O}$ Time Complexity]
Given two functions $f(N)$ and $g(N)$, we say that $f(N)$ is $\mathcal{O}(g(N))$ if there exist constants $c > 0$ and $N_0 \geq 0$ such that $f(N) \leq c \cdot g(N)$ for all $N \geq N_0$.
\end{definition}

It is well known that the fast Fourier transform (FFT) has a time complexity of $\mathcal{O}(N \log(N))$ in 1D \cite{sauer}. Assuming $d$-dimensional inputs of size $N$, i.e. $N^d$ values, this extends to:
\begin{flalign}
&\mathcal{O}((N^d) \log(N^d)) = \mathcal{O}((N^d) \cdot d \log(N)) = \mathcal{O}((N^d) \log(N)) \notag &
\end{flalign}

\indent All computations for this research are performed in MATLAB, which uses the FFTW library for FFT, DCT, and DST function implementation, which offers fast algorithms that handle each with $\mathcal{O}(N^d\log(N))$ for the d-dimensional inputs mentioned \cite{FJ}. The Fourier Spectral method transforms the spatial domain problem into the Fourier frequency domain, where the linear operator, specifically the Laplacian, becomes diagonal. Due to this transformation, the exponential of the linear operator matrix translates to element-wise operations. Specifically, the matrix exponential \( e^{-kA} \) becomes \( (e^{-k\lambda_1}, e^{-k\lambda_2}, \ldots, e^{-k\lambda_N}) \), where \( \lambda_i \) are the eigenvalues of the linear operator matrix \( A \). Thus, instead of performing computationally expensive full matrix multiplications, we perform element-wise exponentiation and multiplication in the frequency domain. This approach enhances computational efficiency, reducing the complexity from $\mathcal{O}(N^3)$ to $\mathcal{O}(N)$ for each dimension, which aids in simulating high-dimensional CNLSEs efficiently. Because the linear operator $A$ operates in frequency space, a FFT or IFFT is required at each iteration to maintain consistency in operation space. In considering the time complexity, we assume the algorithms are acting on a d-dimensional input with each dimension being size $N$, i.e. $N^d$ values

\indent The initialization step involves computing the initial Pad$\acute{\mathrm{e}}$ matrices needed and computing initial FFTs. The complexity of this step is mainly dominated by the FFT operations, which have a complexity of $\mathcal{O}(N^d \log N)$ for FFTs on a grid of $N^d$ values. The main loop of both the Krogstad-P22 and IFRK4-P13 algorithms includes:
\begin{itemize}
    \item \textbf{Computing FFTs and IFFTs:} Each FFT operation in $d$-dimensions has a complexity of $\mathcal{O}((N^d) \log N)$.
    \item \textbf{Matrix exponentials and multiplications:} Matrix operations are performed element-wise, involving $\mathcal{O}(N^d)$ operations per step.
    \item \textbf{Nonlinear function evaluations:} Each function evaluation involves $\mathcal{O}(N^d)$ operations.
\end{itemize}
Combining these, the complexity per iteration of the main loop is $\mathcal{O}((N^d) \log N)$ for both methods. The final computation involves an inverse FFT to transform the solution back to the physical space, which also has a complexity of $\mathcal{O}((N^d) \log N)$. Thus, the overall runtime complexity for both the Krogstad-P22 and IFRK4-P13 algorithms is $\mathcal{O}((N^d) \log N)$.

\subsection{Stability}
\begin{figure}[H]
\begin{minipage}[b]{0.49\linewidth} 
\centering
\centerline{\includegraphics[width=\linewidth,height=\textheight,keepaspectratio]{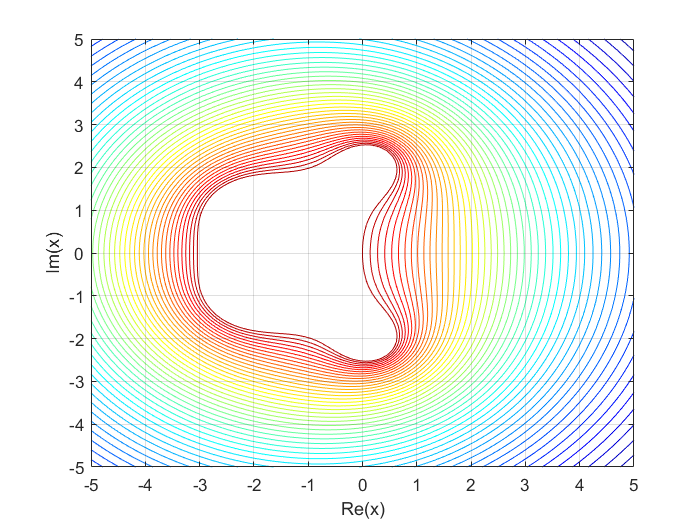}}
\centerline{\text{(a)}}
\end{minipage}
\begin{minipage}[b]{0.49\linewidth} 
\centering
\centerline{\includegraphics[width=\linewidth,height=\textheight,keepaspectratio]{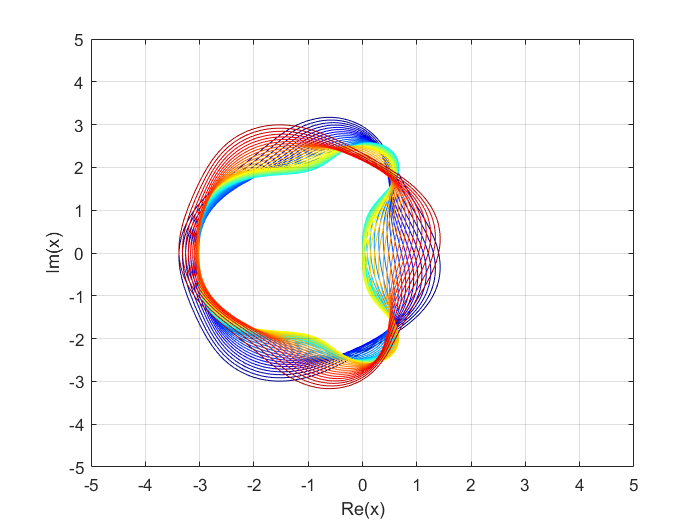}}
\centerline{\text{(b)}}
\end{minipage}
\caption{\footnotesize{Stability regions for equispaced (a) $y\in \R^-$ and (b) $y\in [-2\pi i,2\pi i]$.}}
\label{fig: stability}
\end{figure}

\indent In this section, we will discuss the linear stability analysis of the proposed methods. The stability analysis for the Krogstad-P22 method has already been conducted, so readers are referred to \cite{HBK} for details. For the IFRK4-P13 method, we will perform the linear stability analysis by plotting stability regions, following the approach outlined in \cite{Cox, NZZ} and the references therein. In analyzing the stability of the IFRK4-P13 method, let us consider the following linear equation: 
\begin{flalign} 
    &\frac{du}{dt}=-cu+\lambda u.&
\label{eq:transfsode}
\end{flalign}
Let $x=\lambda k$ and $y=-ck$, where k is the time step. Applying the IFRK4-P13 to Equation \eqref{eq:transfsode} yields the following amplification factor of the method:
\begin{flalign}
    &\frac{u_{n+1}}{u_n}=r(x,y)=c_0+c_1x+c_2x^2+c_3x^3+c_4x^4,
    \label{eq:amp}&
\end{flalign}

where
\begin{flalign}
    &c_0 = 1+y+\frac{1}{2}y^2+\frac{1}{6}y^3+\frac{1}{24}y^4+\frac{1}{96}y^5+O(y^6) \notag&\\
    &c_1 = \frac{17}{24}+\frac{71}{96}y+\frac{587}{1536}y^2+\frac{1205}{9216}y^3+\frac{9847}{294912}y^4+\frac{2965}{393216}y^5+O(y^6) \notag&\\
    &c_2 = \frac{17}{48}+\frac{71}{192}y+\frac{587}{3072}y^2+\frac{1205}{18432}y^3+\frac{9847}{589824}y^4+\frac{903}{262144}y^5+O(y^6) &\\
    &c_3 = \frac{1}{6}+\frac{1}{6}y+\frac{1}{12}y^2+\frac{1}{36}y^3+\frac{1}{144}y^4+\frac{13}{9216}y^5+O(y^6) \notag&\\
    &c_4 = \frac{1}{24}+\frac{1}{24}y+\frac{1}{48}y^2+\frac{1}{144}y^3+\frac{1}{576}y^4+\frac{13}{36864}y^5+O(y^6) \notag&
\end{flalign}
\indent We obtain the boundaries of the stability regions by substituting $u_n=e^{in\theta}$, $\theta \in [0,2\pi]$, which yields $r=e^{i \theta}$, into Equation \eqref{eq:amp} and solving for $x$. However, since we do not know the explicit formula for $|r(x,y)|=1$, we can only plot it, as follows in Figure \ref{fig: stability}.

\indent In Figure \ref{fig: stability} (a), we plotted the stability regions of the IFRK4-P13 for different values of $y \in [-5,0]$ in the complex $x-$ plane. As is discussed in \cite{HBK}, it is known that the stability region must grow as $y$ approaches $-\infty$ for the method to be useful. We can see from Figure \ref{fig: stability} (a) that the stability regions for IFRK4-P13 grow larger as $y$ approaches $-\infty$ while maintaining the shape, suggesting the method is stable.
\\\indent In Figure \ref{fig: stability} (b), we see the stability region of the IFRK4-P13 method for the case where $\lambda$ is complex and $c$ is purely imaginary. We can see from Figure \ref{fig: stability} (b) that the stability region of the method almost includes an interval of imaginary values. This observation suggests an indication of the stability of the IFRK4-P13 method to solve nonlinear dispersive PDEs.

\section{Numerical experiments and discussions}
\label{sec:4}
All following tests are performed using MATLAB 2023b on a Lenovo ThinkPad desktop with i7-10510U CPU and 16 GB memory.
\subsection{Conserved Quantities}
In this section, we will utilize the $L_\infty$ norm error to measure the accuracy of the methods in comparison. We also consider the conservation of mass and energy of simulated nonlinear waves. Following previous related discussion and derivations, \cite{KLS, WXY, IA, BK} we define mass and energy conservation as follows:
\\Mass conservation \cite{IA}:
\begin{flalign}
    &I_i = \int_{-\infty}^{\infty} |\Psi_i|^{2} \, dx&
\label{eq:mass}
\end{flalign}

\noindent Energy conservation \cite{WXY}:
\begin{flalign}
     &E(t) = \frac{1}{2\mu}\int_{R} \sum_{i=1}^{M}{\overline{\Psi_i}(-\Delta)\Psi_i} \, dx - \frac{1}{4}\int_{R} (b\sum_{i=1}^{M}{|\Psi_i|^{4}})+2e(\sum_{i=1}^{M/2}{|\Psi_{2i-1}|^{2}|\Psi_{2i}|^{2}}) \, dx,&\\
    &b=\sigma_{jm} \text{ for } j\neq m, e=\sigma_{jm} \text{ for } j=m. \notag&
     \label{eq:energy}
\end{flalign}
To measure the rate of convergence of the methods we compare, we utilized the following formula:
\begin{flalign}
    &\text{Order}=\frac{\log(||\Psi_{j_n}-\Psi_{j_{n-1}}||_\infty / ||\Psi_{j_{n+1}}-\Psi_{j_n}||_\infty)}{\log(2)}&
\end{flalign}

In our numerical experiments, we will simulate the solution profile of 1D, 2D, and 3D M-CNLSEs under various boundary conditions. These experiments aim to compare the effectiveness in terms of accuracy and computational efficiency of the Krogstad-P22 method with the IFRK4-P13 method for approximating CNLSEs.

\subsection{1D CNLSE}
\noindent {\bf Example 1 (Single soliton propagation):}\ 
\\In this first instance, we consider a CNLSE problem with a known analytical solution as a benchmark problem to analyze and compare the accuracy and efficiency of the proposed methods. We consider system \eqref{eq:1} for $M=2$ to model single soliton propagation on the domain $\Omega=[-20,80]$, represented as follows:
\begin{flalign}
&i\frac{\partial \Psi_1}{\partial t}+\frac{1}{\mu}\frac{\partial^2 \Psi_1}{\partial x^2} +(\sigma_{11}|\Psi_1|^2 + \sigma_{12}|\Psi_2|^2)\Psi_1 = 0,\quad x_L<x<x_R \label{eq:m2singsys}&&
\\&i\frac{\partial \Psi_2}{\partial t}+\frac{1}{\mu}\frac{\partial^2 \Psi_2}{\partial x^2} +(\sigma_{21}|\Psi_1|^2 + \sigma_{22}|\Psi_2|^2)\Psi_2 = 0,\quad x_L<x<x_R \notag&&
\end{flalign}
together with periodic boundary conditions and the following initial conditions:
\begin{flalign}
&\Psi_1(x,0) = \sqrt{\frac{\mu \alpha}{1+e}}\text{sech}(\sqrt{\mu\alpha}x)\exp{(ivx)}, \label{eq:singinit}&&
\\&\Psi_2(x,0) = \sqrt{\frac{\mu \alpha}{1+e}}\text{sech}(\sqrt{\mu\alpha}x)\exp{(ivx)},\notag&&
\end{flalign}
where $\alpha_1$,$\alpha_2$=$\frac{1}{\mu}$, $\sigma_{jj}=1$ for $j=1,2$, $\sigma_{jm}=e$ for $j\neq m$, and $\alpha$, $\mu$, $e$, and $v$ being constants. Considering the derivation by Wadati et al. in \cite{WIH}, the exact solution for this system is known to be:
\begin{flalign}
\Psi_j (x,t) = \sqrt{\frac{\mu \alpha}{1+e}}\text{sech}(\sqrt{\mu\alpha}(x-vt))\exp{(i(vx-[\textstyle \frac{v^2}{2}-\alpha]t))}, \text{ for }j=1,2.&&
\label{eq:m2singsysexact}
\end{flalign}

In this example, we wish to confirm empirically that the Krogstad-P22 and IFRK4-P13 methods are fourth-order convergent in time. To test this convergence, we found the maximum norm error by comparing the numerical solution with the given exact solution for decreasing time steps starting with $k=0.025$, with the following parameters corresponding to linearly birefringent fibers:
$N=1024$, $T=5$, $v=1$, $\alpha=1$, $\mu=2$, and $e$=2/3. Our choice of $N=1024$ is to avoid any errors coming from the spatial discretization, and the computed results are given in Table \ref{tab:Table 1}. 

From Table \ref{tab:Table 1}, we can see that both methods converge to fourth-order accuracy as expected. However, we can notice that Krogstad-P22 acquires errors within $10^{-7}$, $10^{-8}$, and $10^{-11}$ in less CPU time than the IFRK4-P13 method does. This suggests Krogstad-22 acquires better accuracy with larger time step $k$ than IFRK4-P13, and suits this problem better.
\begin{table}[H]
\caption{Comparison of temporal convergence for Krogstad-P22 and IFRK4-P13 methods} 
\begin{center} \footnotesize
\resizebox{\linewidth}{!}{
\begin{tabular}{ccccc|cccc}
 \hline
  &\multicolumn{3}{c}{Krogstad-P22}& &\multicolumn{3}{c}{IFRK4-P13}\\
  \\\cmidrule{1-4}\cmidrule{6-8}
 $k$ & $\norm{\Psi_1}_\infty$ & Order &CPU-time && $\norm{\Psi_1}_\infty$ & Order &CPU-time\\
  \hline\\
$\frac{1}{40}$& 3.9012$\times 10^{-7}$ & - & 0.122 && 7.0262$\times 10^{-6}$ & - & 0.146\\\\
$\frac{1}{80}$& 1.8594$\times 10^{-8}$ &4.3910 &0.226 && 3.3418$\times 10^{-7}$ & 4.3940 & 0.241\\\\
$\frac{1}{160}$& 9.8323$\times 10^{-10}$ &4.2411 &0.494 && 1.9149$\times 10^{-8}$ & 4.1253 & 0.442\\\\
$\frac{1}{320}$& 5.5945$\times 10^{-11}$ &4.1354 &0.817 && 1.2039$\times 10^{-9}$ & 3.9916 & 0.798\\\\
$\frac{1}{640}$& 3.3592$\times 10^{-12}$  &4.0578 &1.352 && 7.5121$\times 10^{-11}$ & 4.0023 & 1.484\\\\
\hline
\end{tabular}}
\end{center} 
\label{tab:Table 1} 
\end{table} 

\begin{figure}[H]
\begin{minipage}[b]{0.33\linewidth} 
\centering
\centerline{\includegraphics[width=\linewidth,height=\textheight,keepaspectratio]{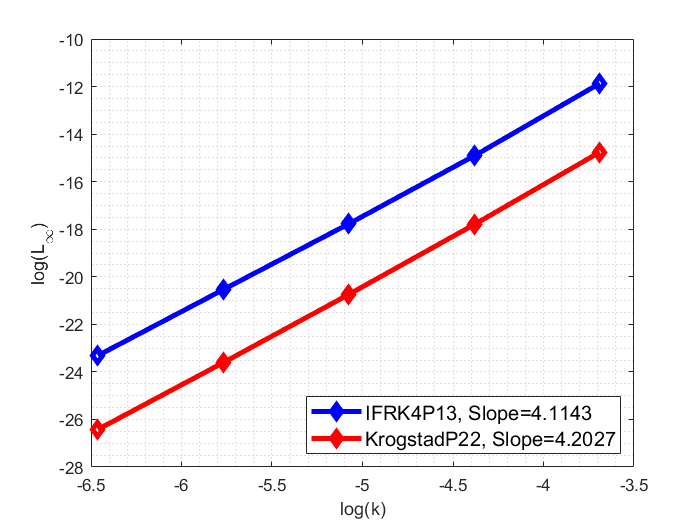}}
\centerline{\text{(a) Periodic Boundary}}
\end{minipage}
\begin{minipage}[b]{0.33\linewidth}
\centering
\centerline{\includegraphics[width=\linewidth,height=\textheight,keepaspectratio]{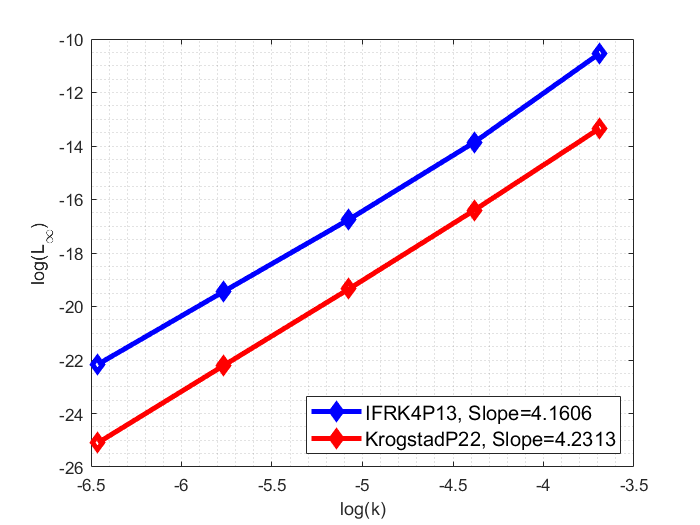}}
\centerline{\text{(b) Neumann Boundary}}
\end{minipage}
\begin{minipage}[b]{0.33\linewidth}
\centering
\centerline{\includegraphics[width=\linewidth,height=\textheight,keepaspectratio]{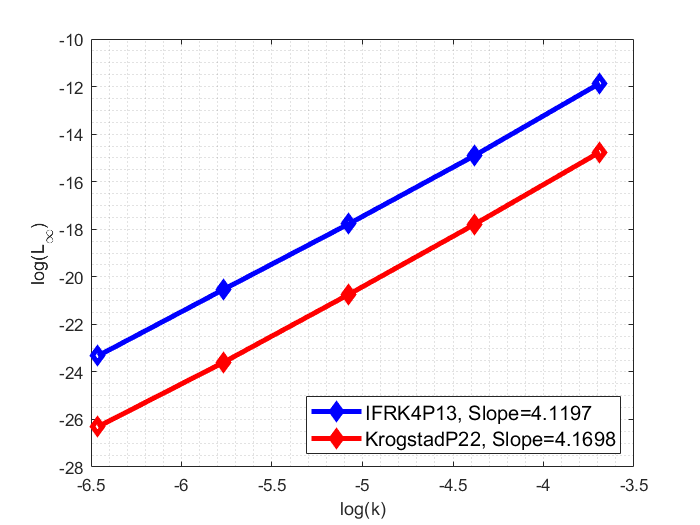}}
\centerline{\text{(c) Dirichlet Boundary}}
\end{minipage}
\caption{\footnotesize{Temporal convergence of the Krogstad-P22 method vs IFRK4-P13 under various boundary conditions}}
\label{fig: fig1}
\end{figure}

The data presented in Table \ref{tab:Table 1}, along with the temporal convergence under homogeneous Dirichlet and Neumann boundary conditions using the same parameters, are shown on a log-log scale in Figure \ref{fig: fig1}. From Figure \ref{fig: fig1}, it is evident that both methods achieved the expected fourth-order convergence in the temporal direction for all the boundary conditions considered.

\begin{figure}[H]
\begin{minipage}[b]{0.33\linewidth} 
\centering
\centerline{\includegraphics[width=\linewidth,height=\textheight,keepaspectratio]{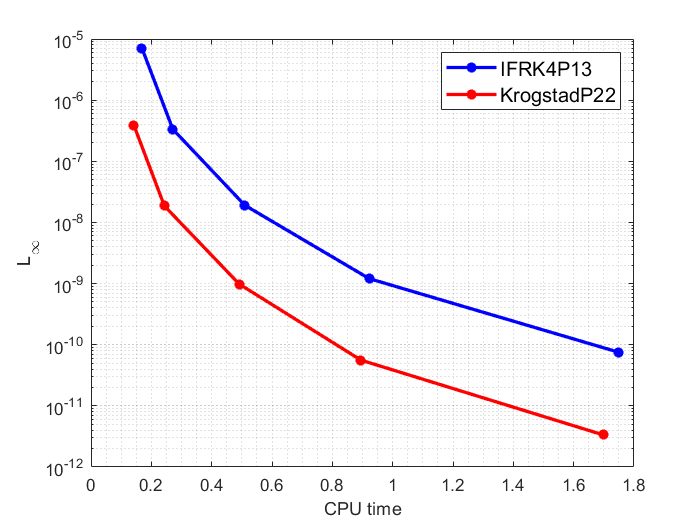}}
\centerline{\text{(a) Periodic Boundary}}
\end{minipage}
\begin{minipage}[b]{0.33\linewidth}
\centering
\centerline{\includegraphics[width=\linewidth,height=\textheight,keepaspectratio]{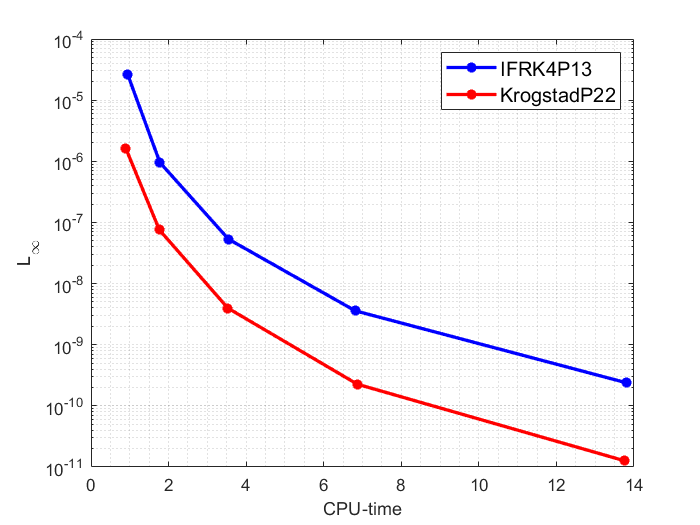}}
\centerline{\text{(b) Neumann Boundary}}
\end{minipage}
\begin{minipage}[b]{0.33\linewidth}
\centering
\centerline{\includegraphics[width=\linewidth,height=\textheight,keepaspectratio]{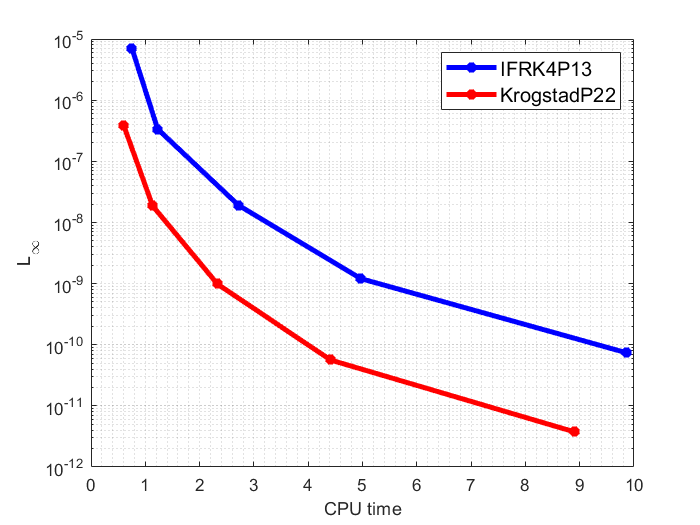}}
\centerline{\text{(c) Dirichlet Boundary}}
\end{minipage}
\caption{\footnotesize{Comparison of the computational efficiency of Krogstad-P22 and IFRK4-P13 methods}}
\label{fig: fig3}
\end{figure}

\noindent In Figure \ref{fig: fig3}, we compared the computational efficiency of the proposed methods. From Figure \ref{fig: fig3}, we observe that the Krogstad-P22 method outperformed the IFRK4-P13 method in terms of efficiency and accuracy.

\begin{figure}[H]
\begin{minipage}[b]{0.33\linewidth} 
\centering
\centerline{\includegraphics[width=\linewidth,height=\textheight,keepaspectratio]{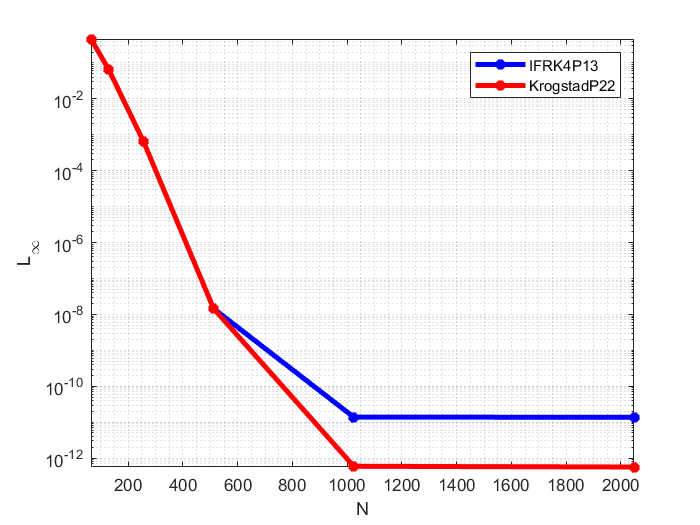}}
\centerline{\text{(a) Periodic Boundary}}
\end{minipage}
\begin{minipage}[b]{0.33\linewidth}
\centering
\centerline{\includegraphics[width=\linewidth,height=\textheight,keepaspectratio]{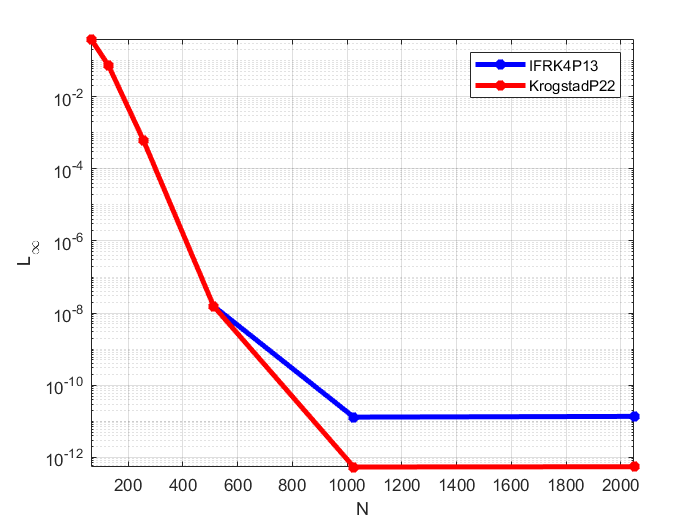}}
\centerline{\text{(b) Neumann Boundary}}
\end{minipage}
\begin{minipage}[b]{0.33\linewidth}
\centering
\centerline{\includegraphics[width=\linewidth,height=\textheight,keepaspectratio]{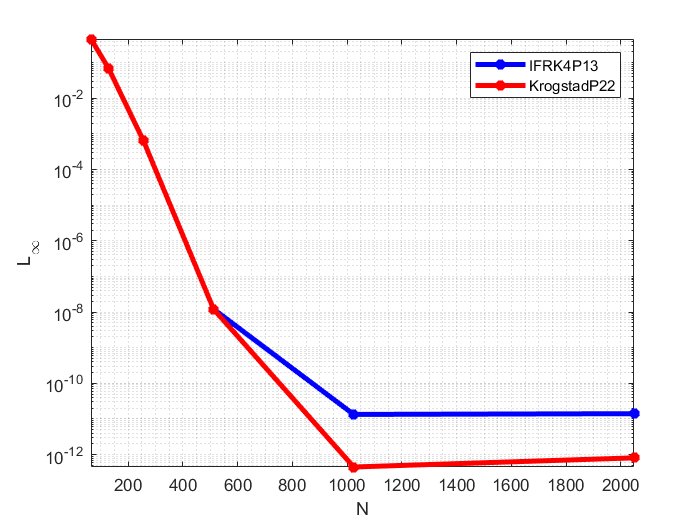}}
\centerline{\text{(c) Dirichlet Boundary}}
\end{minipage}
\caption{\footnotesize{Spatial convergence of the Fourier spectral method in combination with Krogstad-P22 and IFRK4-P13 under various boundary conditions}}
\label{fig: fig2}
\end{figure}

We conducted another set of experiments to determine whether the Fourier spectral method achieves spectral convergence in the spatial direction under periodic, homogeneous Dirichlet, and Neumann boundary conditions. In these experiments, we ran the simulations until $T = 1$, keeping $k = 0.001$ fixed and repeatedly doubling $N = 64$ for each simulation, storing the maximum errors. The results are displayed in Figure \ref{fig: fig2}. From Figure \ref{fig: fig2}, it is evident that the Fourier spectral method exhibited the expected spectral convergence in the spatial direction for each boundary condition.

We also considered the accuracy of both methods under the periodic boundary condition by comparing the mass conservation of soliton polarization $|\Psi_1|$ and system energy conservation, with parameters $x_L=-20, x_R=80$ $N=1024$, $k=0.0125$, $v=1$, $\alpha=1$, $\mu=2$, and $e=2/3$ up to $T=40$. 

\begin{figure}[H]
\begin{minipage}[b]{0.33\linewidth} 
\centering
\centerline{\includegraphics[width=\linewidth,height=\textheight,keepaspectratio]{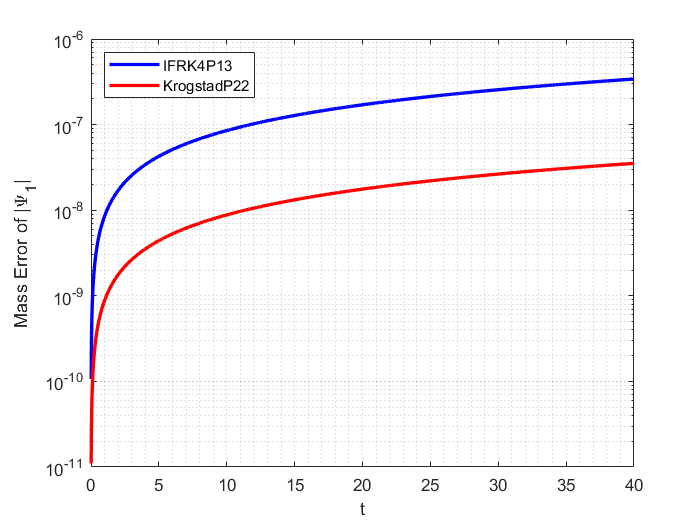}}
\centerline{\text{(a) Periodic Boundary}}
\end{minipage}
\begin{minipage}[b]{0.33\linewidth}
\centering
\centerline{\includegraphics[width=\linewidth,height=\textheight,keepaspectratio]{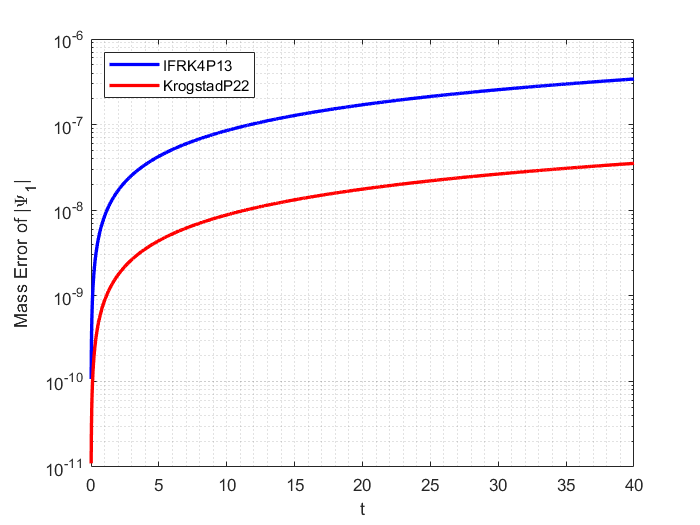}}
\centerline{\text{(b) Neumann Boundary}}
\end{minipage}
\begin{minipage}[b]{0.33\linewidth}
\centering
\centerline{\includegraphics[width=\linewidth,height=\textheight,keepaspectratio]{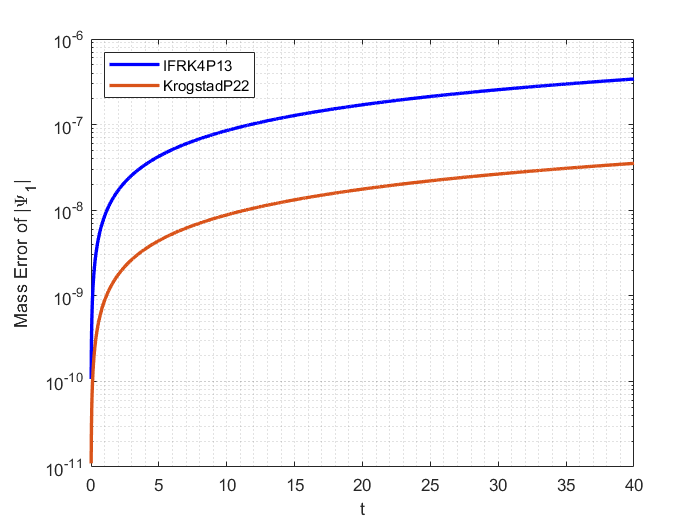}}
\centerline{\text{(c) Dirichlet Boundary}}
\end{minipage}
\caption{\footnotesize{Mass conservation comparison under various boundary conditions.}}
\label{fig: fig4}
\end{figure}

In Figure \ref{fig: fig4}, we compared the absolute error of the mass of $|\Psi_1|$ for the Krogstad-P22 and IFRK4-P13 methods, using Equation \eqref{eq:mass}, under various boundary conditions. We can see that the Krogstad-P22 method is conserving mass better over a longer time than the IFRK4-P13 method for each boundary condition. 

\begin{figure}[H]
\begin{minipage}[b]{0.33\linewidth} 
\centering
\centerline{\includegraphics[width=\linewidth,height=\textheight,keepaspectratio]{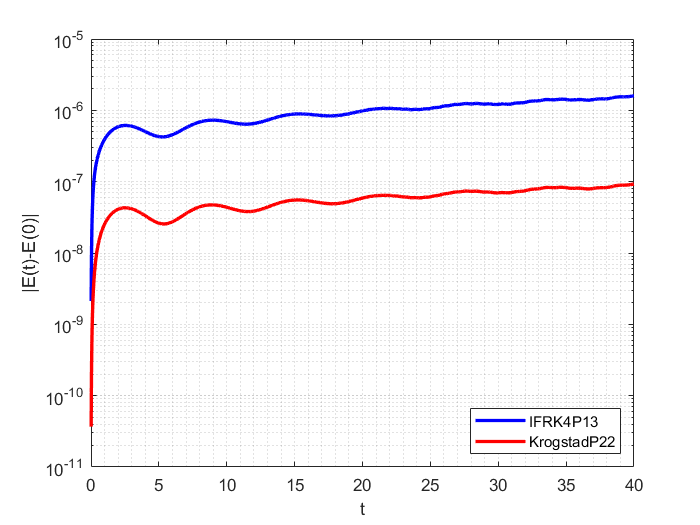}}
\centerline{\text{(a) Periodic Boundary}}
\end{minipage}
\begin{minipage}[b]{0.33\linewidth}
\centering
\centerline{\includegraphics[width=\linewidth,height=\textheight,keepaspectratio]{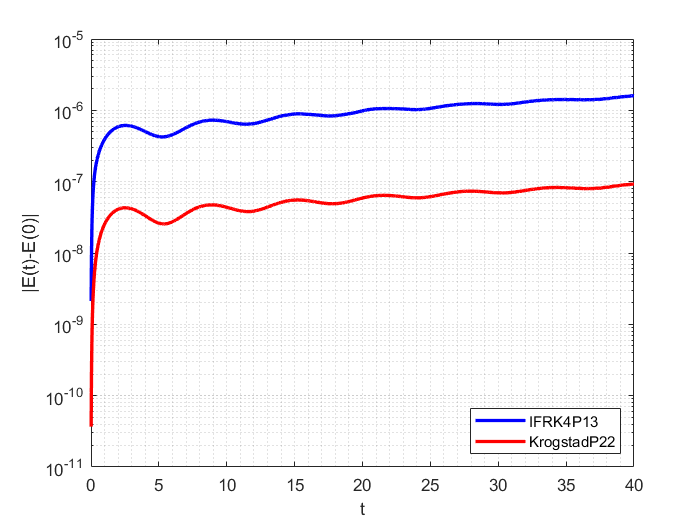}}
\centerline{\text{(b) Neumann Boundary}}
\end{minipage}
\begin{minipage}[b]{0.33\linewidth}
\centering
\centerline{\includegraphics[width=\linewidth,height=\textheight,keepaspectratio]{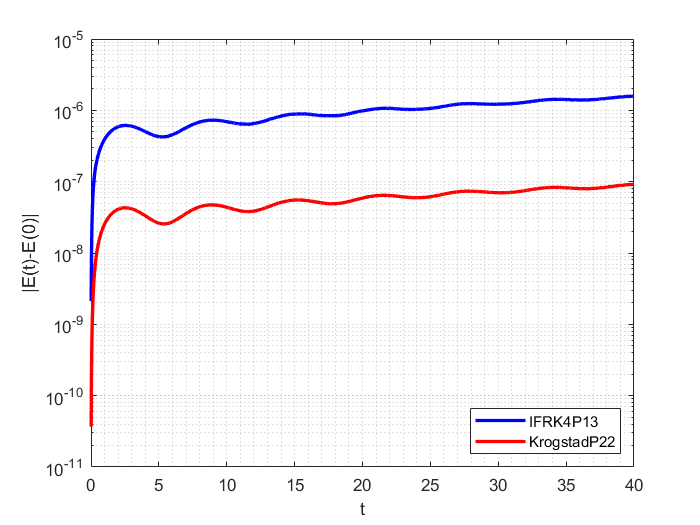}}
\centerline{\text{(c) Dirichlet Boundary}}
\end{minipage}
\caption{\footnotesize{Energy conservation comparison under various boundary conditions.}}
\label{fig: fig5}
\end{figure}

\noindent Figure \ref{fig: fig5} compares the absolute error of the energy of the whole soliton system, computed using Equation \eqref{eq:energy}, for the same methods under the same set of boundary conditions as in Figure \ref{fig: fig4}. We can see that Krogstad-P22 conserves system energy better over a longer time than the IFRK4-P13 method. From these results, we can infer that while neither method conservative, they both conserve mass and energy well over a longer time with Krogstad-P22 conserving these properties better.

\noindent {\bf Example 2 (Interaction of two solitons):} 
\\\indent In this instance, we analyze a CNLSE problem modeling the interaction of two solitons to evaluate the methods' ability to preserve expected interaction behavior. We solve system \eqref{eq:m2singsys} on the domain $\Omega=[-40,40]$ with $\mu=1$, simulating two solitons of different amplitudes moving in opposite directions. The initial conditions, taken from \cite{MIT, Yang}, are as follows:
\begin{flalign}
&\Psi_1(x,0) = \sqrt{2}r_1 \text{sech}(r_1x+x_{10})\exp{(iv_1x)} ,\notag&&
\\&\Psi_2(x,0) = \sqrt{2}r_2 \text{sech}(r_2x+x_{20})\exp{(iv_2x)},
\label{eq:2intinit}
\end{flalign}
with homogeneous Neumann boundary conditions:
\begin{flalign}
\frac{\partial \Psi_1(x,t)}{\partial x}=\frac{\partial \Psi_2(x,t)}{\partial x}=0, \text{ at $x=x_L,x_R$ for $t\geq0$},&&
\label{eq:hnbc}
\end{flalign}
where $\sigma_{jj}$=1, $\sigma_{12}=\sigma_{21}=e$, $r_1$=1.2, $r_2$=1, $v_1=-v_2=\frac{v}{4}$. For all computations in this example, we use $x_L$=-40, and $x_R$=40. Note that our choices of $x_{10}$ and $x_{20}$ are somewhat arbitrary. Given that $x_{10}>$0 and $x_{20}<$0 are large enough, they should have no effect on the results of our collision \cite{BK,Yang}.

In table \ref{tab:Table 3}, we confirmed that the Krogstad-P22 and IFRK4P13 methods remain fourth-order convergent in time under the homogeneous Neumann boundary condition, using $x_{10}=30$, $x_{20}=-10$, and halving time steps starting with $k=0.0125$, where $N=512$, $e$=1, and $v$=1. Our choices for $e$ and $v$ are motivated by their influence on the CNLSEs we consider. $e$=1 and $v$=1 yields Manakov's equations, which are fully integrable, and provide an elastic interaction \cite{Yang}. Our choice of $N = 512$ works nicely for the Fourier spectral method, requiring less computations compared to a larger choice of N in exchange for accuracy.
\begin{table}[H]
\begin{center} \footnotesize
\resizebox{\linewidth}{!}{
\begin{tabular}{cccc|cccc}
 \hline
  &\multicolumn{3}{c}{Krogstad-P22}&&\multicolumn{3}{c}{IFRK4P13}
  \\\cmidrule{2-4}\cmidrule{6-8}
 $k$ & $\norm{\Psi_1}_\infty$ & Order &CPU-time && $\norm{\Psi_1}_\infty$ & Order &CPU-time\\
  \hline\\
  $\frac{1}{160}$& 1.0610$\times 10^{-7}$ & - & 1.946 && 1.2244$\times 10^{-6}$ & - & 2.120\\\\
  $\frac{1}{320}$& 4.8585$\times 10^{-9}$ & 4.4488 & 3.772 && 7.1846$\times 10^{-8}$ & 4.0910 & 3.923\\\\
  $\frac{1}{640}$& 2.9108$\times 10^{-10}$ & 4.0610 & 7.545 && 4.2224$\times 10^{-9}$ & 4.0888 & 10.128\\\\
  $\frac{1}{1280}$& 2.0024$\times 10^{-12}$ & 3.8616 & 15.170 && 2.5493$\times 10^{-10}$ & 4.0247 & 22.391\\\\
\hline
\end{tabular}}
\end{center} 
\caption{Comparison of temporal convergence for Krogstad-P22 and IFRK4P13 methods}
\label{tab:Table 3} 
\end{table} 
From Table \ref{tab:Table 3}, we can see that both methods converge to fourth-order accuracy as expected. However, we can also notice that Krogstad-P22 acquires errors within $10^{-9}$ and $10^{-10}$ in less CPU time than the IFRK4P13 method does. In general, Krogstad-P22 acquires better accuracy sooner than IFRK4P13.

\begin{figure}[H]
\begin{minipage}[b]{0.33\linewidth} 
\centering
\centerline{\includegraphics[width=\linewidth,height=\textheight,keepaspectratio]{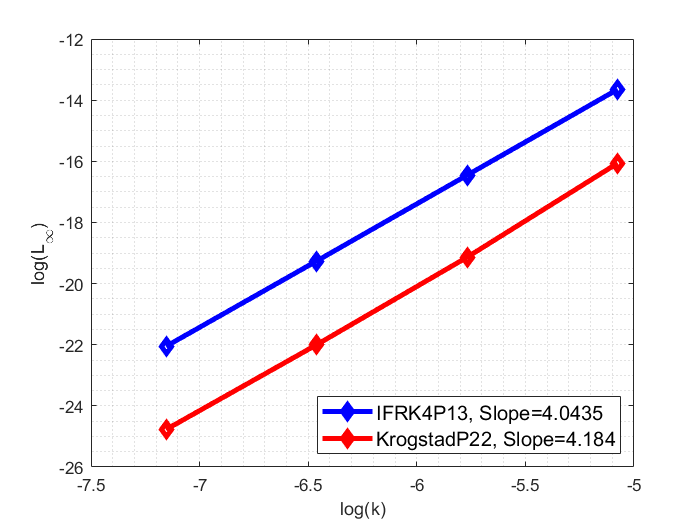}}
\centerline{\text{(a) Periodic Boundary}}
\end{minipage}
\begin{minipage}[b]{0.33\linewidth}
\centering
\centerline{\includegraphics[width=\linewidth,height=\textheight,keepaspectratio]{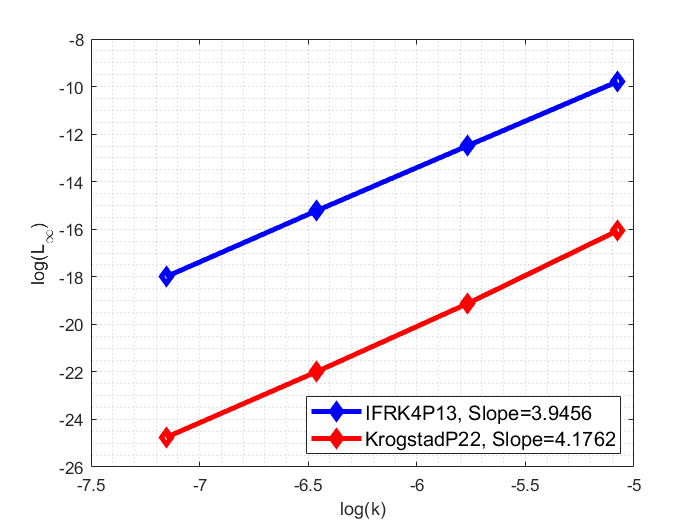}}
\centerline{\text{(b) Neumann Boundary}}
\end{minipage}
\begin{minipage}[b]{0.33\linewidth}
\centering
\centerline{\includegraphics[width=\linewidth,height=\textheight,keepaspectratio]{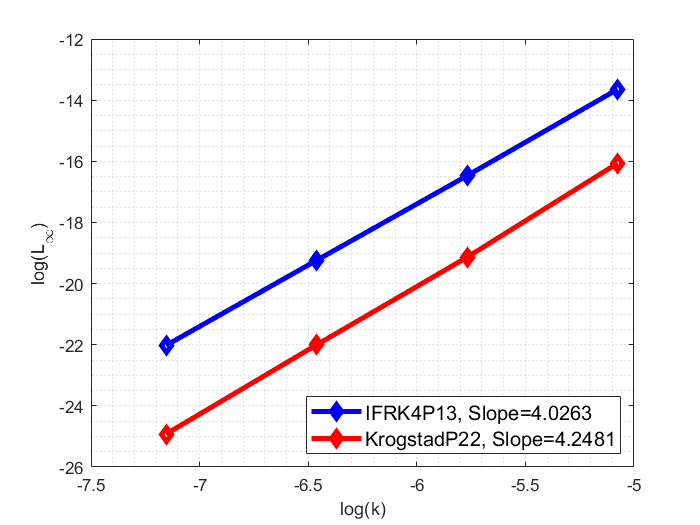}}
\centerline{\text{(c) Dirichlet Boundary}}
\end{minipage}
\caption{\footnotesize{Temporal convergence of the Krogstad-P22 method vs IFRK4P13 under various boundary conditions}}
\label{fig: fig6}
\end{figure}

The data presented in Table \ref{tab:Table 3}, along with the temporal convergence under periodic and homogeneous Dirichlet boundary conditions using the same parameters, are shown on a log-log scale in Figure \ref{fig: fig6}. From Figure \ref{fig: fig6}, it is evident that both methods continue to exhibit the expected fourth-order convergence in the temporal direction for all the boundary conditions considered.

\begin{figure}[H]
\begin{minipage}[b]{0.33\linewidth} 
\centering
\centerline{\includegraphics[width=\linewidth,height=\textheight,keepaspectratio]{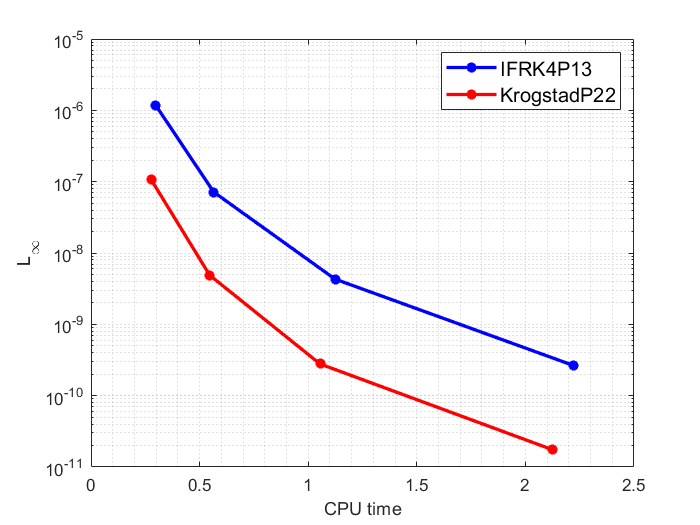}}
\centerline{\text{(a) Periodic Boundary}}
\end{minipage}
\begin{minipage}[b]{0.33\linewidth}
\centering
\centerline{\includegraphics[width=\linewidth,height=\textheight,keepaspectratio]{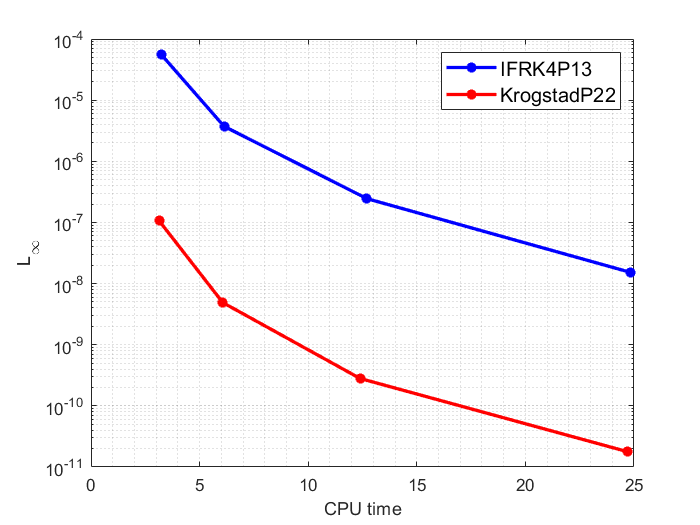}}
\centerline{\text{(b) Neumann Boundary}}
\end{minipage}
\begin{minipage}[b]{0.33\linewidth}
\centering
\centerline{\includegraphics[width=\linewidth,height=\textheight,keepaspectratio]{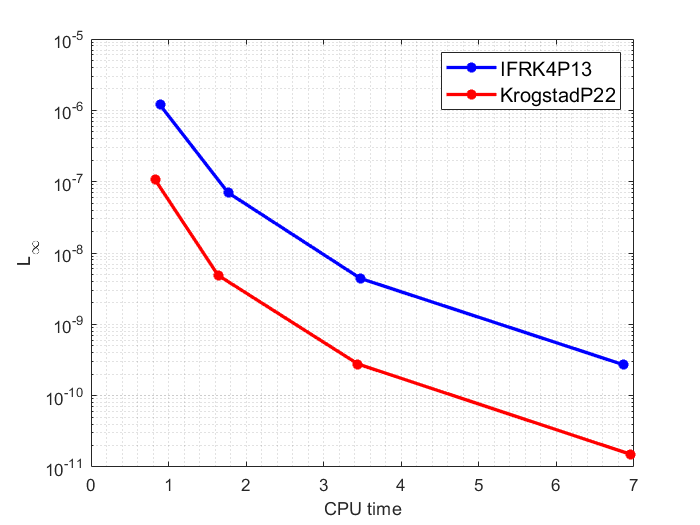}}
\centerline{\text{(c) Dirichlet Boundary}}
\end{minipage}
\caption{\footnotesize{CPU time and $L_\infty$ error of $|\Psi_1|$ for decreasing $k$ under various boundary conditions}}
\label{fig: fig7}
\end{figure}

 In Figure \ref{fig: fig7}, we compared the computational efficiency of the proposed methods using the data from Table \ref{tab:Table 3}. From Figure \ref{fig: fig7}, we observe that the Krogstad-P22 method outperformed the IFRK4P13 method in terms of efficiency and accuracy, which agrees with our results from Table \ref{tab:Table 3}.
\begin{table}[H]
\setlength{\tabcolsep}{30pt}
\centering
\begin{tabular}{cccc}
\hline 
T & $I_1$ & $I_2$ & $E(t)$ \\ \hline 
0 & 4.800000 & 4.0000000 & 4.196860  \\
10 & 4.800000 & 4.0000000 & 4.196812  \\
20 & 4.800000 & 4.0000000 & 4.196871  \\
30 & 4.800000 & 4.0000000 & 4.196880  \\
40 & 4.800000 & 4.0000000 & 4.196837  \\
50 & 4.800000 & 4.0000000 & 4.201194  \\
60 & 4.800000 & 4.0000000 & 4.196873  \\
70 & 4.800000 & 4.0000000 & 4.196866  \\ 
80 & 4.800000 & 4.0000000 & 4.196864  \\ 
90 & 4.800000 & 4.0000000 & 4.196843  \\ 
100 & 4.800000 & 4.0000000 & 4.196892 \\ 
\hline
\end{tabular}
\caption{Mass and Energy of $\Psi_1$ and $\Psi_2$}
\label{tab: Table 4}
\end{table}
\indent We conducted another simulation using the Krogstad-P22 method under the initial conditions of this example, and analyze the conserved quantities, with $e=1$, $v=1$, $x_{10}=-x_{20}=30$, $N=1024$, and $k=0.01$ for $T=100$. In Table \ref{tab: Table 4}, we consider the conserved quantities. From the table, we observed the mass of the first soliton $|\Psi_1|$ is $I_1\approx$4.800000 and the mass of the second soliton $|\Psi_2|$ is $I_1\approx$4.0000000 up to $T=100$. We also observed that the energy of the soliton system, given by $E(t)$ in Table \ref{tab: Table 4}, fluctuates during soliton collision, but returns to an accuracy within four decimal places by final $T$=100. The conserved quantities in Table \ref{tab: Table 4} demonstrate that the Krogstad-P22 method conserves soliton system mass and energy well over a longer time. 

\begin{figure}[H]
\begin{minipage}[b]{0.5\linewidth} 
\centering
\centerline{\includegraphics[width=\linewidth,height=\textheight,keepaspectratio]{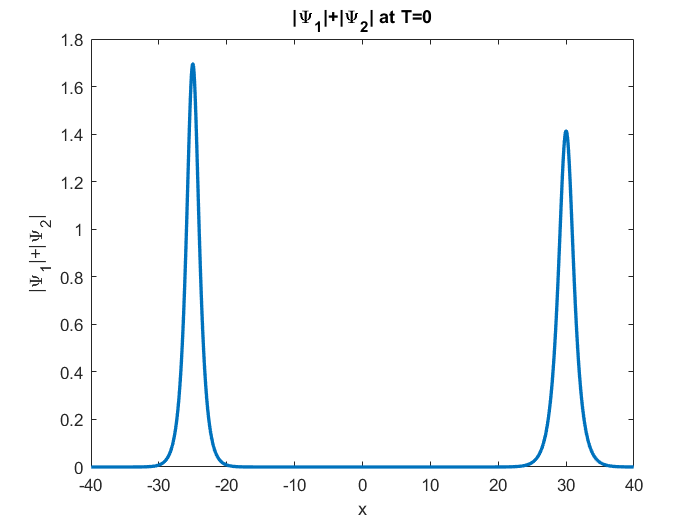}}
\end{minipage}
\begin{minipage}[b]{0.5\linewidth}
\centering
\centerline{\includegraphics[width=\linewidth,height=\textheight,keepaspectratio]{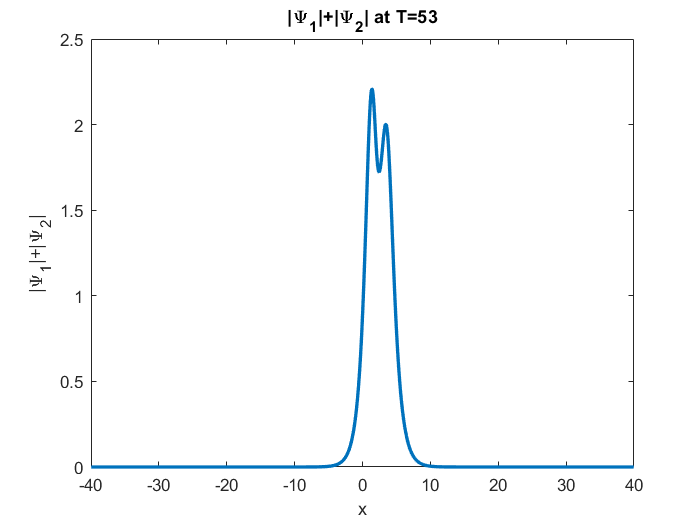}}
\end{minipage}
\begin{minipage}[b]{0.5\linewidth}
\centering
\centerline{\includegraphics[width=\linewidth,height=\textheight,keepaspectratio]{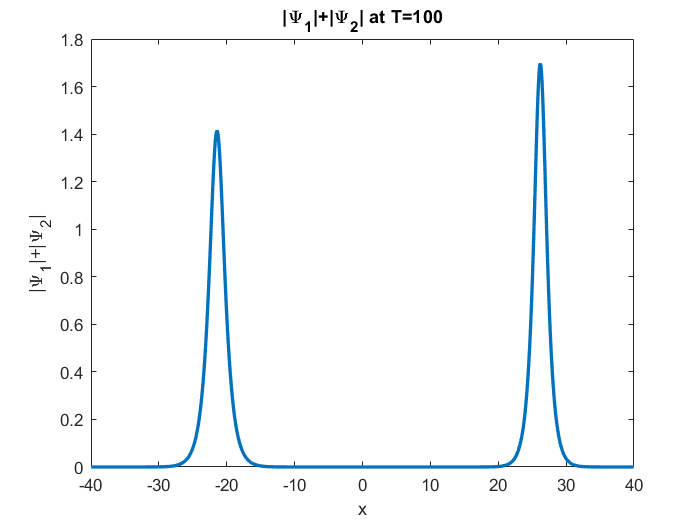}}
\end{minipage}
\begin{minipage}[b]{0.5\linewidth}
\centering
\centerline{\includegraphics[width=\linewidth,height=\textheight,keepaspectratio]{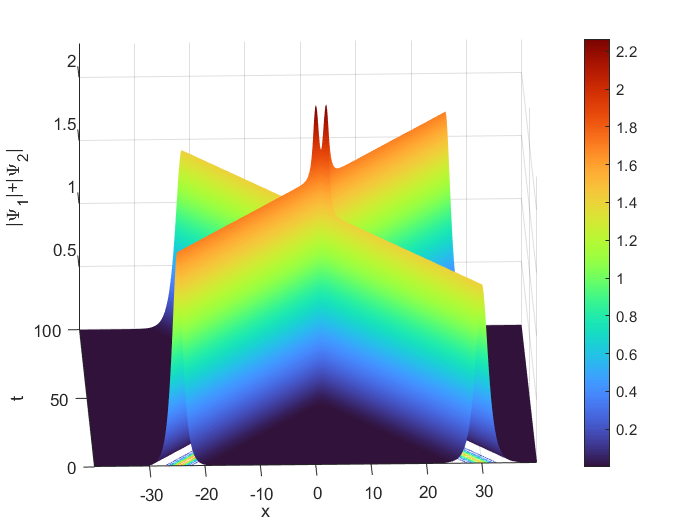}}
\end{minipage}
\caption{\footnotesize{Numerical simulation of elastic soliton collision modeled from system \eqref{eq:m2singsys} with conditions \eqref{eq:2intinit} }}
\label{fig: fig10}
\end{figure}

Figure \ref{fig: fig10} illustrates the collision scenario detailed in Table \ref{tab: Table 4}. In Figure \ref{fig: fig10}, the simulation is shown at times before, during, and after the collision under elastic parameters. We observed the amplitudes of the solitons increase briefly when colliding and then return to their previous amplitudes afterward, continuing their movements as before the collision. This soliton interaction is elastic, and agrees with the analytical predictions in \cite{Yang}.

\begin{figure}[H]
    \centering
    \begin{minipage}[b]{0.49\textwidth}
        \centering
        \begin{subfigure}[b]{\textwidth}
            \includegraphics[width=\textwidth]{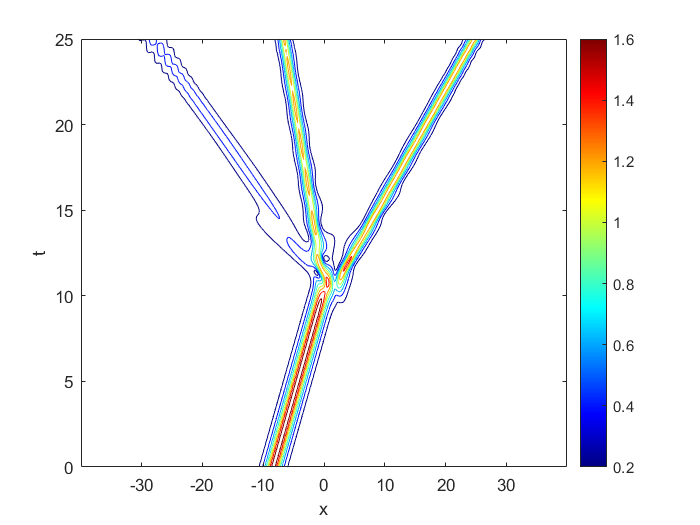}
            \caption{}
        \end{subfigure}

        \begin{subfigure}[b]{\textwidth}
            \includegraphics[width=\textwidth]{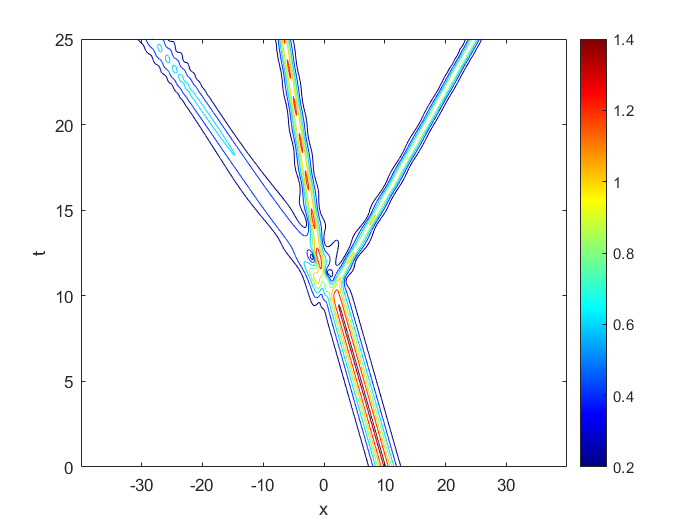}
            \caption{}
        \end{subfigure}
        \caption{Solitons $|\Psi_1|$ (a) and $|\Psi_2|$ (b) colliding, with $e=3$ and $v=1.6$}
        \label{fig: fig11}
    \end{minipage}
    \hfill 
    \begin{minipage}[b]{0.49\textwidth}
        \centering
        \begin{subfigure}[b]{\textwidth}
            \includegraphics[width=\textwidth]{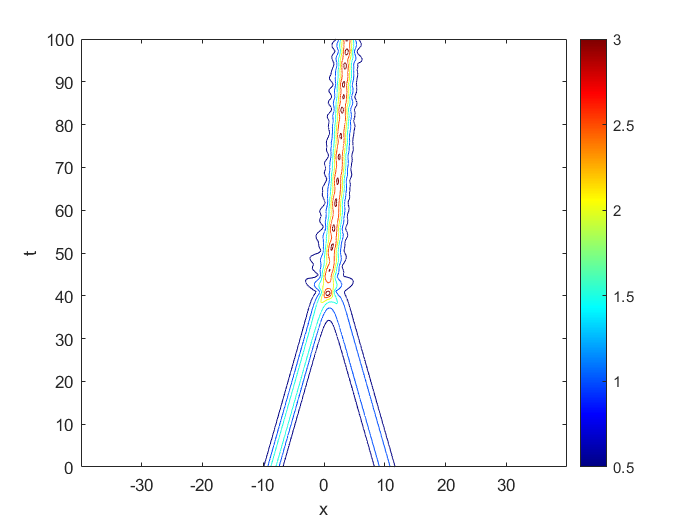}
            \caption{}
        \end{subfigure}

        \begin{subfigure}[b]{\textwidth}
            \includegraphics[width=\textwidth]{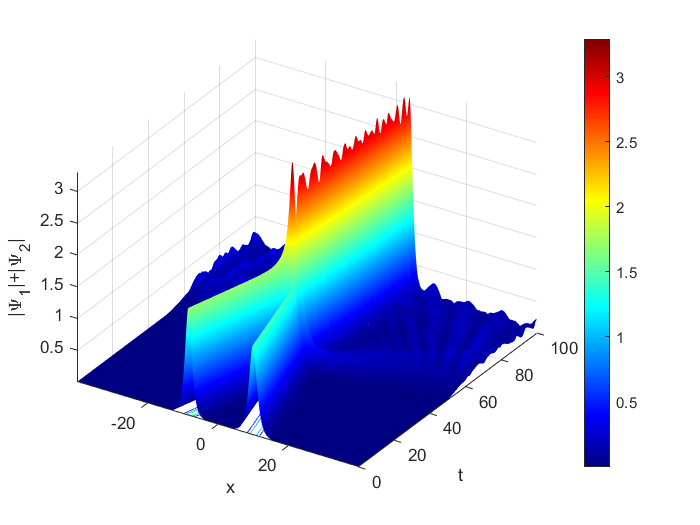}
            \caption{}
        \end{subfigure}
        \caption{Contour plot (a) and surface plot (b) of solitons $|\Psi_1|+|\Psi_2|$ colliding, with $e=0.3$ and $v=0.4$}
        \label{fig: fig12}
    \end{minipage}
\end{figure}

To observe that the Krogstad-P22 method models the interaction of two solitons as expected, we simulate two inelastic collisions discussed in \cite{Yang}. Both simulations have $N=800$ and $k=0.01$. For the first inelastic simulation, we used $x_{20}=-x_{10}=10$ with $e=3$ and $v=1.6$. The contour plots of the first soliton $|\Psi_1|$ and the second soliton $|\Psi_2|$ are shown in Figure \ref{fig: fig11} (a) and (b), respectively. We see in Figure \ref{fig: fig11} that the soliton collision yields the creation of a new soliton, which can happen when $e$ is small with slow velocity $v$. For the second inelastic simulation, we used the same parameters except with $e=0.3$ and $v=0.4$, modeled in Figure \ref{fig: fig12}. We see in Figure \ref{fig: fig12} that the solitons merge into a single soliton upon collision, which can happen when $e$ is large and velocity $v$ is moderate. These inelastic interactions agree nicely with the analytical predictions for these test parameters discussed by Yang in \cite{Yang}. The results from Table \ref{tab: Table 4}, and Figures \ref{fig: fig10}, \ref{fig: fig11}, and \ref{fig: fig12} suggest that the Krogstad-P22 method accurately approximates soliton solutions.

\newpage
\noindent {\bf Example 3 (Interaction of four solitons):} \\
\indent In this instance, we want to analyze a CNLSE problem modeling the interaction of four solitons in a nonlinear medium. We consider the CNLSE on the domain $\Omega=[-40,40]$ represented by system \eqref{eq:1} for M=4 as in \cite{BK}, with $\alpha_1,\alpha_2=1$, $\sigma_{jm}=b$ for $j=m$, and $\sigma_{jm}=e$ for $j \neq m$, as follows:
\begin{flalign}
&i\frac{\partial \Psi_1}{\partial t}+\frac{\partial^2 \Psi_1}{\partial x^2} +(b|\Psi_1|^2 + e(|\Psi_2|^2+|\Psi_3|^2+|\Psi_4^2))\Psi_1 = 0,\quad \notag&&
\\&i\frac{\partial \Psi_2}{\partial t}+\frac{\partial^2 \Psi_2}{\partial x^2} +(b|\Psi_2|^2 + e(|\Psi_1|^2+|\Psi_3|^2+|\Psi_4^2))\Psi_2 = 0,\quad&&
\\&i\frac{\partial \Psi_3}{\partial t}+\frac{\partial^2 \Psi_3}{\partial x^2} +(b|\Psi_3|^2 + e(|\Psi_1|^2+|\Psi_2|^2+|\Psi_4^2))\Psi_3 = 0,\quad\notag&&
\\&i\frac{\partial \Psi_4}{\partial t}+\frac{\partial^2 \Psi_4}{\partial x^2} +(b|\Psi_4|^2 + e(|\Psi_1|^2+|\Psi_3|^2+|\Psi_3^2))\Psi_4 = 0,\quad\notag&&
\\ &x_L<x<x_R\notag
\label{eq:int4eq}
\end{flalign}
together with initial conditions:
\begin{flalign}
&\Psi_1(x,0)=\sqrt{2}r_1\text{sech}(r_1x+x_{10})\exp{(iv_1x)},\notag& \\ 
&\Psi_2(x,0)=\sqrt{2}r_2\text{sech}(r_2x-x_{10})\exp{(-iv_2x)},\notag& \\ 
&\Psi_3(x,0)=\sqrt{2}r_3\text{sech}(r_3x+x_{30})\exp{(iv_3x)},& \\ 
&\Psi_4(x,0)=\sqrt{2}r_4\text{sech}(r_4x-x_{30})\exp{(-iv_4x)},\notag& 
\label{eq:int4eqinit}
\end{flalign}
and homogeneous Dirichlet boundary conditions:
\begin{flalign}
\Psi_1(x,t)=\Psi_2(x,t)=0, \text{ at $x=x_L,x_R$ for $t\geq0$},&&
\end{flalign}
where $x_{10}=10$, $x_{30}=30$, $b=1$, $e=1$, $r_1=1.0$, $r_2=1.2$, $r_3=1.3$, $r_4=1.4$, $v_1=v_2=\frac{1}{8}$, and $~v_3=v_4=\frac{1}{4}$. 

We performed the temporal convergence of the Krogstad-P22 and IFRK4-P13 methods for decreasing time step k, starting with $k=0.025$, with our results shown in Table \ref{tab:Table 5}. From this Table, we can see that both methods converge to fourth-order in time as expected. We can also notice that Krogstad-P22 acquires errors within $10^{-9}$, and $10^{-10}$ in less CPU time than the IFRK4-P13 method does, and generally acquires better accuracy sooner than IFRK4-P13 in this case.

\begin{table}[H]
\begin{center} \footnotesize
\resizebox{\linewidth}{!}{
\begin{tabular}{ccccc|ccc}
 \hline
  &\multicolumn{3}{c}{Krogstad-P22}&&\multicolumn{3}{c}{IFRK4-P13}
  \\\cmidrule{2-4}\cmidrule{5-8}
 $k$ & $\norm{\Psi_1}_\infty$ & Order &CPU-time && $\norm{\Psi_1}_\infty$ & Order &CPU-time\\
  \hline\\
  $\frac{1}{80}$& 1.8917$\times 10^{-7}$ & - & 1.357 && 2.7016$\times 10^{-6}$ & - & 1.397 \\\\
  $\frac{1}{160}$& 1.0115$\times 10^{-8}$ & 4.2251 & 2.678 && 1.8358$\times 10^{-7}$ & 3.879 & 2.835 \\\\
  $\frac{1}{320}$& 5.9274$\times 10^{-10}$ & 4.0930 & 5.001 && 1.1699$\times 10^{-8}$ & 3.9720 & 5.021 \\\\
  $\frac{1}{640}$& 3.5027$\times 10^{-11}$ & 4.0808 & 13.817 && 7.2855$\times 10^{-10}$ & 4.0052 & 15.261 \\\\
\hline
\end{tabular}}
\end{center} 
\caption{Comparison of temporal convergence for Krogstad-P22 and IFRK4-P13 methods in interaction of four solitons}
\label{tab:Table 5} 
\end{table} 

\begin{figure}[H]
\begin{minipage}[b]{0.33\linewidth} 
\centering
\centerline{\includegraphics[width=\linewidth,height=\textheight,keepaspectratio]{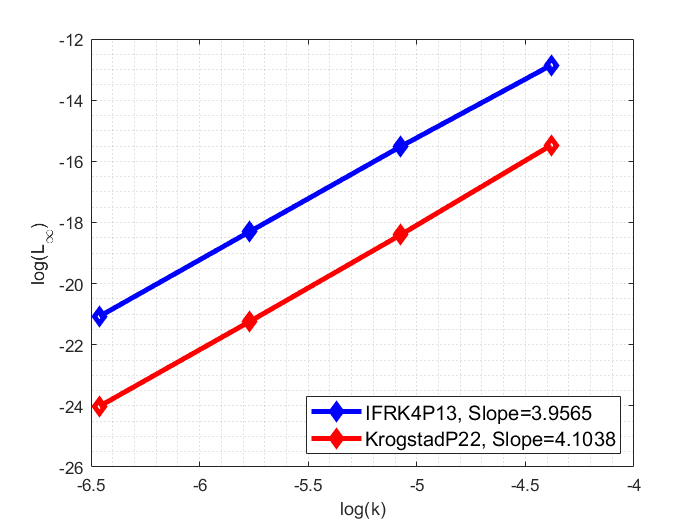}}
\centerline{\text{(a) Periodic Boundary}}
\end{minipage}
\begin{minipage}[b]{0.33\linewidth}
\centering
\centerline{\includegraphics[width=\linewidth,height=\textheight,keepaspectratio]{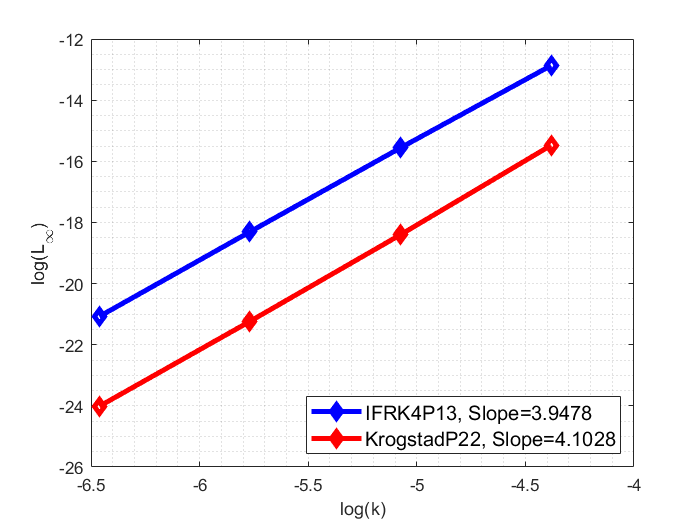}}
\centerline{\text{(b) Neumann Boundary}}
\end{minipage}
\begin{minipage}[b]{0.33\linewidth}
\centering
\centerline{\includegraphics[width=\linewidth,height=\textheight,keepaspectratio]{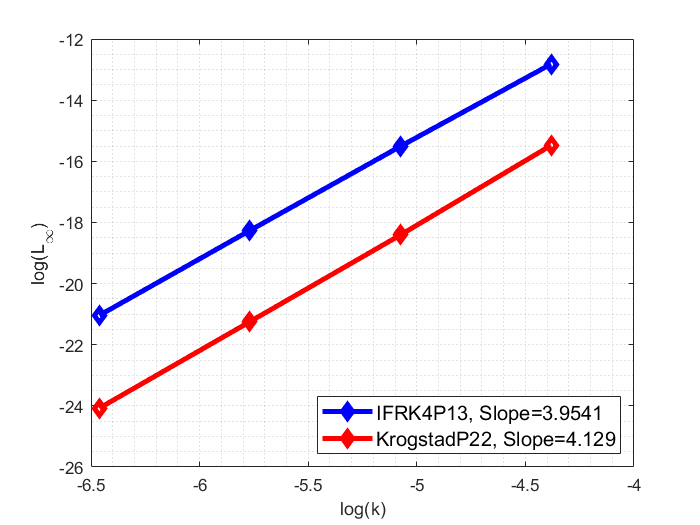}}
\centerline{\text{(c) Dirichlet Boundary}}
\end{minipage}
\caption{\footnotesize{Temporal convergence of the Krogstad-P22 method vs IFRK4-P13 under various boundary conditions}}
\label{fig: fig13}
\end{figure}

The data presented in Table \ref{tab:Table 5}, along with the temporal convergence under periodic and homogeneous Neumann boundary conditions using the same parameters, are shown on a log-log scale in Figure \ref{fig: fig13}. From Figure \ref{fig: fig13}, it is evident that both methods achieved the expected fourth-order convergence in the temporal direction for all the boundary conditions considered.

\begin{figure}[H]
\begin{minipage}[b]{0.33\linewidth} 
\centering
\centerline{\includegraphics[width=\linewidth,height=\textheight,keepaspectratio]{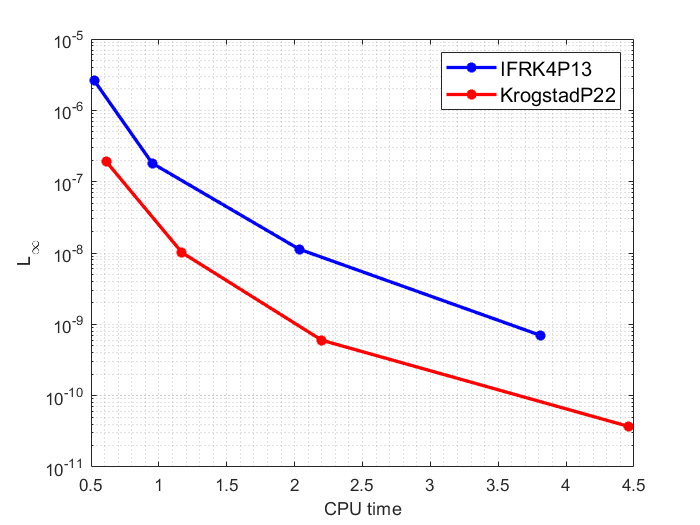}}
\centerline{\text{(a) Periodic Boundary}}
\end{minipage}
\begin{minipage}[b]{0.33\linewidth}
\centering
\centerline{\includegraphics[width=\linewidth,height=\textheight,keepaspectratio]{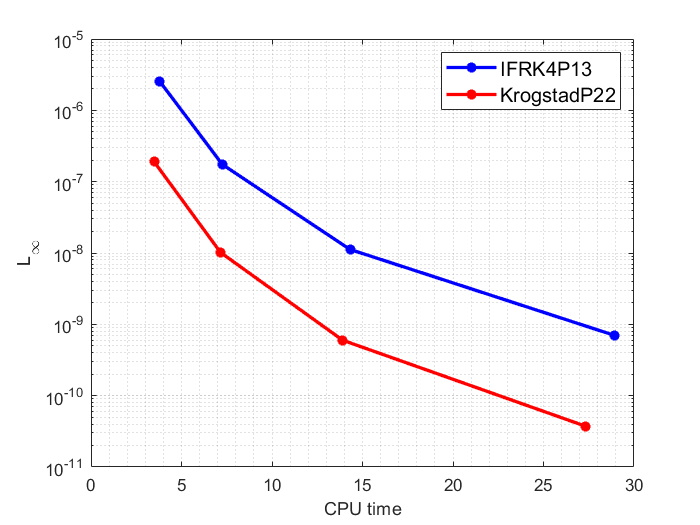}}
\centerline{\text{(b) Neumann Boundary}}
\end{minipage}
\begin{minipage}[b]{0.33\linewidth}
\centering
\centerline{\includegraphics[width=\linewidth,height=\textheight,keepaspectratio]{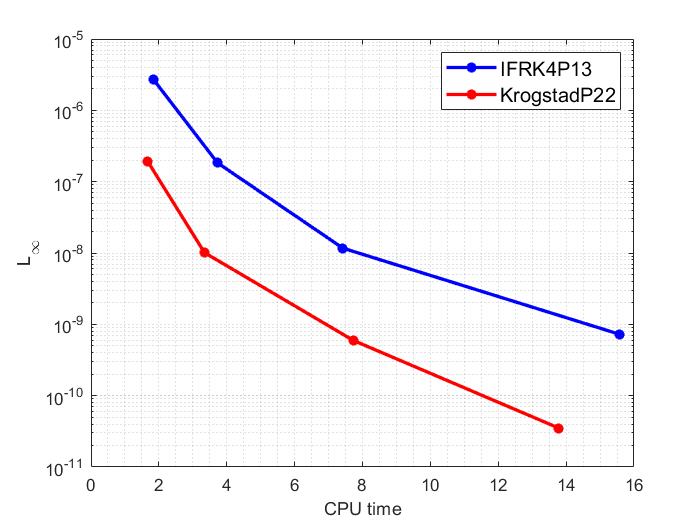}}
\centerline{\text{(c) Dirichlet Boundary}}
\end{minipage}
\caption{\footnotesize{Comparison of computational efficiency of Krogstad-P22 and IFRK4-P13 under various boundary conditions for Example 2.}}
\label{fig: fig14}
\end{figure}

In Figure \ref{fig: fig14}, we compared the $L_\infty$ error of $|\Psi_1|$ against the CPU time required for the simulation as the time step $k$ decreases. Under each boundary condition, we observed that while the Krogstad-P22 method achieves better accuracy than the IFRK4-P13 method at each time step, it also requires more time. Conversely, the IFRK4-P13 method only matches the accuracy of Krogstad-P22 at the next smaller time step, thus taking longer to achieve a similar level of accuracy.

We also considered the energy conservation of system \eqref{eq:int4eq} and mass conservation of $|\Psi_1|$ under the given homogeneous Dirichlet boundary conditions, with parameters $N=800$, $k=0.01$, $b=1$, $e=1$, and $v=1$ up to $T=100$. In Table \ref{tab: Table 6}, the mass of each soliton and the system energy are given under the Krogstad-P22 method. We can see that the mass for each soliton conserves up to at least five decimal places. Although the energy shows greater variation, especially near the soliton collision where the energy briefly decreases, the system energy returns to accuracy up to four decimal places by T=80.
\begin{table}[H]
\setlength{\tabcolsep}{10pt}
\centering
\begin{tabular}{cccccc}
\hline 
T & $I_1$ & $I_2$ & $I_3$ & $I_4$ & $E(t)$ \\
\hline
0 & 4.0 & 4.8 & 5.2 & 5.60000 & 11.6196335 \\
10 & 4.0000000 & 4.800000 & 5.200000 & 5.60000 & 11.6196327 \\
20 & 4.0000000 & 4.800000 & 5.200000 & 5.60000 & 11.6196225 \\
30 & 4.0000000 & 4.800000 & 5.200000 & 5.60000 & 11.6779830 \\
40 & 4.0000000 & 4.800000 & 5.200000 & 5.60000 & 11.9710313 \\
50 & 4.0000000 & 4.800000 & 5.200000 & 5.60000 & 8.6316253 \\
60 & 4.0000000 & 4.800000 & 5.200000 & 5.60000 & 11.5990733 \\
70 & 4.0000000 & 4.800000 & 5.200000 & 5.60000 & 11.6195369 \\
80 & 4.0000000 & 4.800000 & 5.200000 & 5.60000 & 11.6196270 \\
90 & 4.0000000 & 4.800000 & 5.200000 & 5.60000 & 11.6196291 \\
100 & 4.0000000 & 4.800000 & 5.200000 & 5.60000 & 11.6196290 \\
\hline
\end{tabular}
\caption{Mass of $\Psi_i$ and system Energy for Example 2}
\label{tab: Table 6}
\end{table}
\indent

\begin{figure}[H]
\begin{minipage}[b]{0.33\linewidth} 
\centering
\centerline{\includegraphics[width=\linewidth,height=\textheight,keepaspectratio]{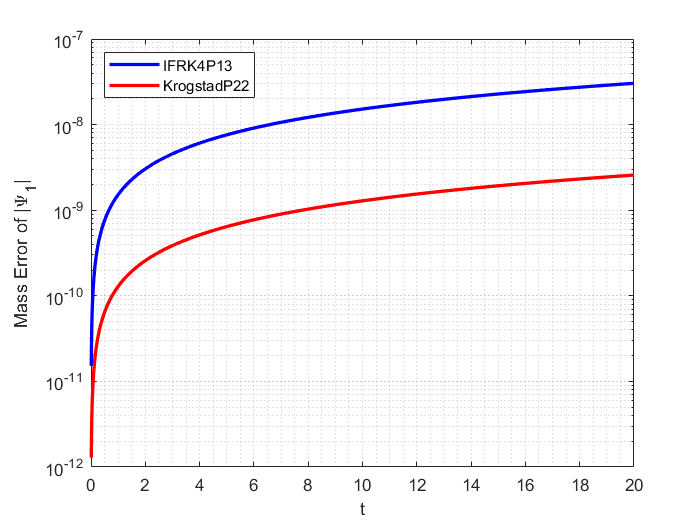}}
\centerline{\text{(a) Periodic Boundary}}
\end{minipage}
\begin{minipage}[b]{0.33\linewidth}
\centering
\centerline{\includegraphics[width=\linewidth,height=\textheight,keepaspectratio]{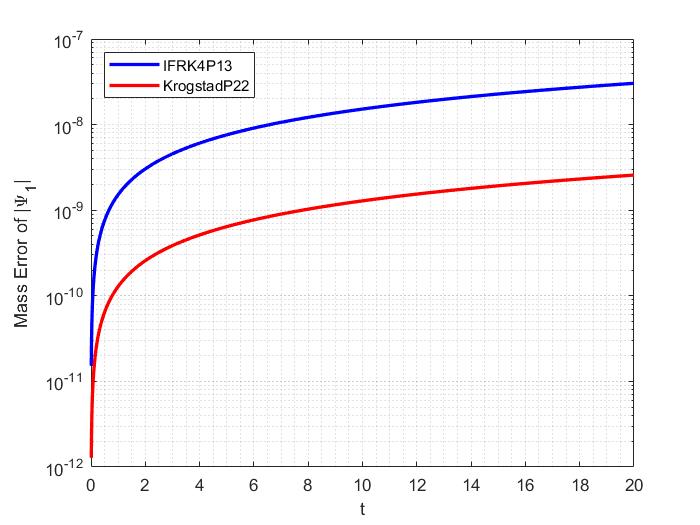}}
\centerline{\text{(b) Neumann Boundary}}
\end{minipage}
\begin{minipage}[b]{0.33\linewidth}
\centering
\centerline{\includegraphics[width=\linewidth,height=\textheight,keepaspectratio]{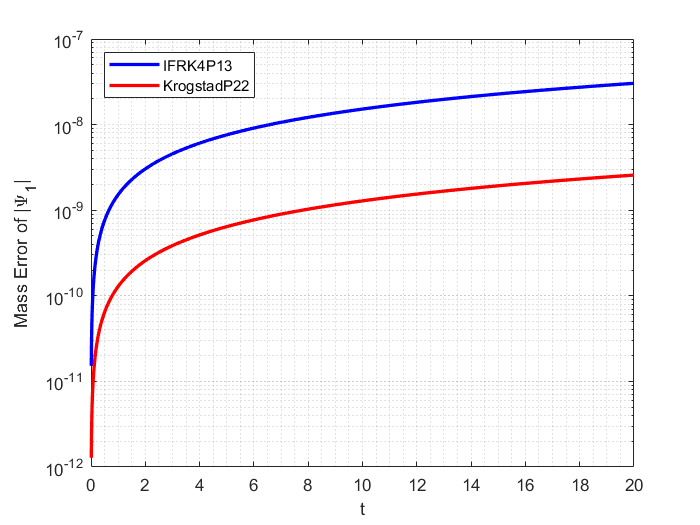}}
\centerline{\text{(c) Dirichlet Boundary}}
\end{minipage}
\caption{\footnotesize{Mass conservation comparison under various boundary conditions.}}
\label{fig: fig15}
\end{figure}

\begin{figure}[H]
\begin{minipage}[b]{0.33\linewidth} 
\centering
\centerline{\includegraphics[width=\linewidth,height=\textheight,keepaspectratio]{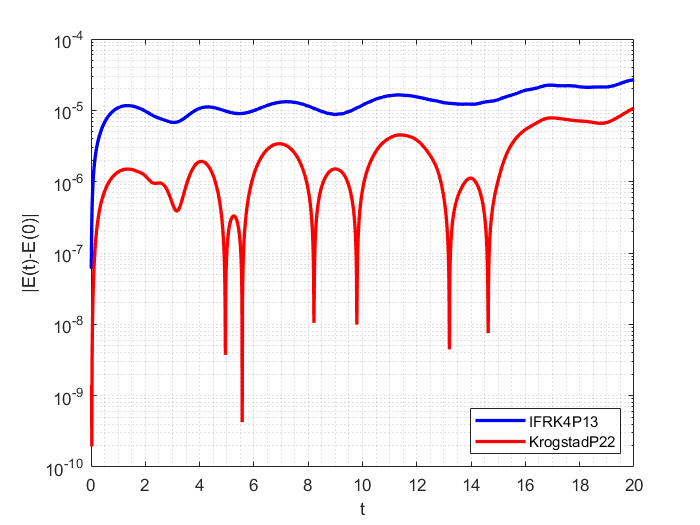}}
\centerline{\text{(a) Periodic Boundary}}
\end{minipage}
\begin{minipage}[b]{0.33\linewidth}
\centering
\centerline{\includegraphics[width=\linewidth,height=\textheight,keepaspectratio]{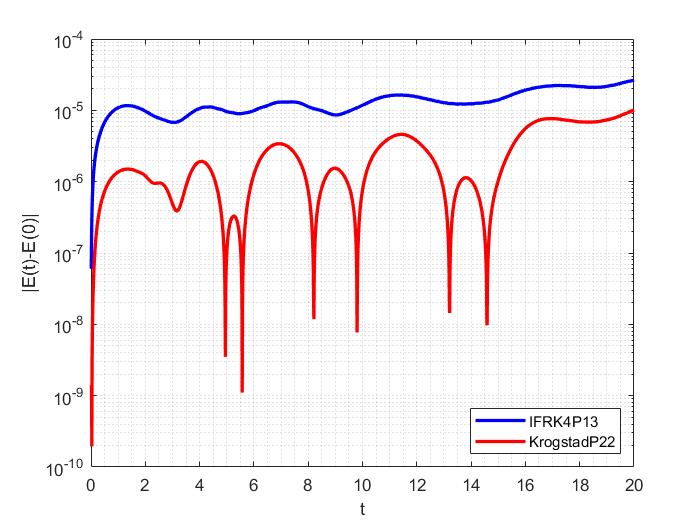}}
\centerline{\text{(b) Neumann Boundary}}
\end{minipage}
\begin{minipage}[b]{0.33\linewidth}
\centering
\centerline{\includegraphics[width=\linewidth,height=\textheight,keepaspectratio]{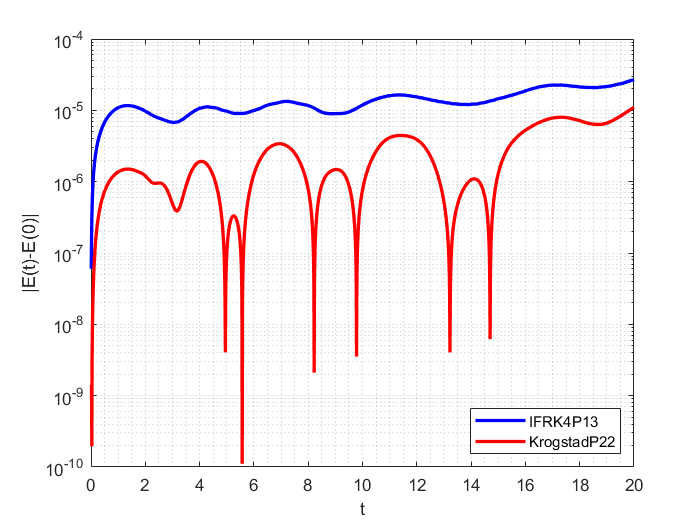}}
\centerline{\text{(c) Dirichlet Boundary}}
\end{minipage}
\caption{\footnotesize{Energy conservation comparison under various boundary conditions.}}
\label{fig: fig16}
\end{figure}

Using the same parameters, the mass and energy conservation under both methods are shown up to $T=20$ in Figure \ref{fig: fig15} and Figure \ref{fig: fig16}, respectively. We choose to show up to $T=20$ for these Figures because the error increases closer to the soliton interactions and the full graphs can make the relative comparison between the two methods difficult to observe.

In Figure \ref{fig: fig15}, we compared the absolute error of the mass of $|\Psi_1|$ under different boundary conditions. We observed that the Krogstad-P22 method conserved mass better over a longer time than the IFRK4-P13 method for each boundary condition. 

In Figure \ref{fig: fig16}, we compared the absolute error of the energy of the whole soliton system under various boundary conditions. From Figure \ref{fig: fig16}, we can see that Krogstad-P22 conserves system energy better over a longer time than the IFRK4-P13 method. From these results, we can infer that while neither method conservative, they both conserved mass and energy well over a longer time with Krogstad-P22 conserving these properties better.

\subsubsection{Blow-up}
\begin{figure}[H]
\begin{minipage}[b]{0.33\linewidth} 
\centering
\centerline{\includegraphics[width=\linewidth,height=\textheight,keepaspectratio]{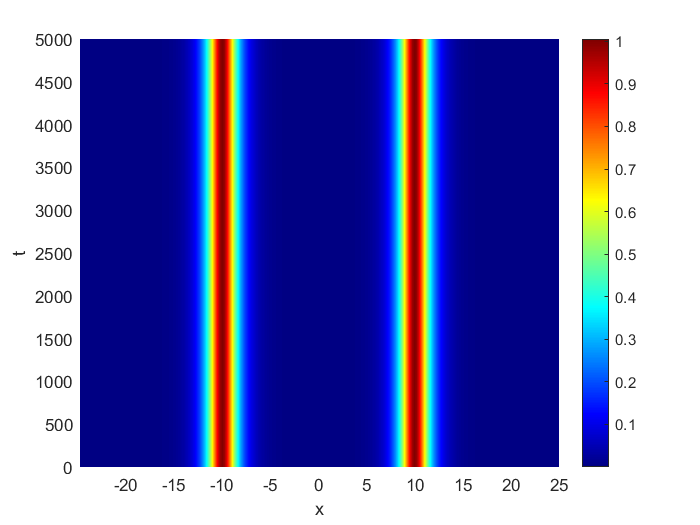}}
\centerline{\text{(a) Periodic Boundary}}
\end{minipage}
\begin{minipage}[b]{0.33\linewidth}
\centering
\centerline{\includegraphics[width=\linewidth,height=\textheight,keepaspectratio]{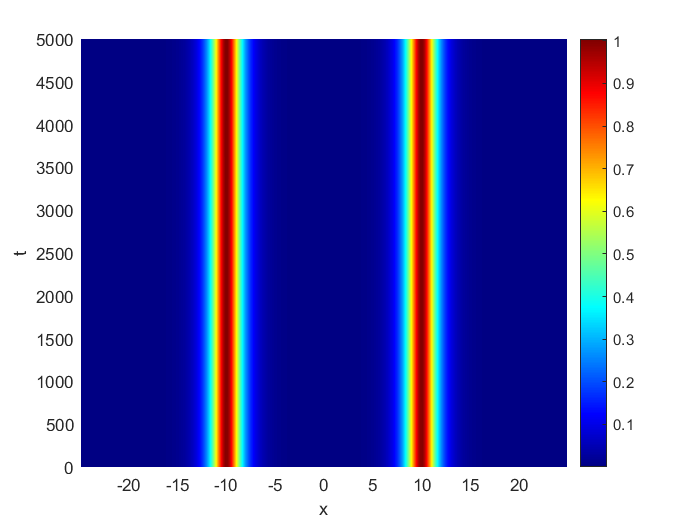}}
\centerline{\text{(b) Neumann Boundary}}
\end{minipage}
\begin{minipage}[b]{0.33\linewidth}
\centering
\centerline{\includegraphics[width=\linewidth,height=\textheight,keepaspectratio]{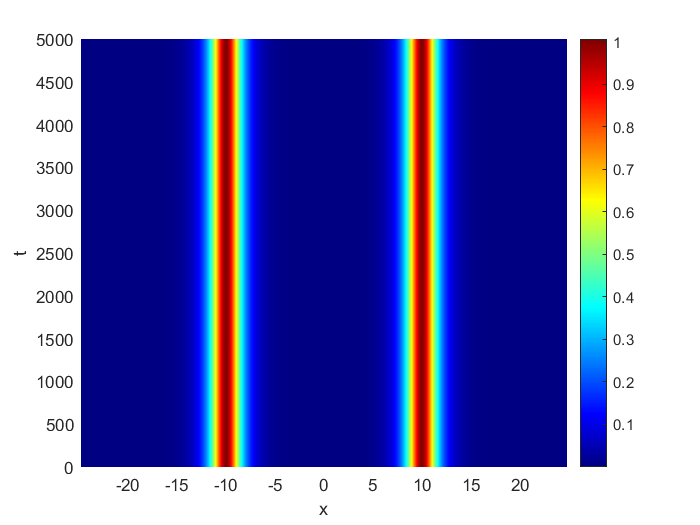}}
\centerline{\text{(c) Dirichlet Boundary}}
\end{minipage}
\caption{\footnotesize{Empirical blow-up analysis of the Krogstad-P22 method}}
\label{fig: fig17}
\end{figure}

\begin{figure}[H]
\begin{minipage}[b]{0.33\linewidth} 
\centering
\centerline{\includegraphics[width=\linewidth,height=\textheight,keepaspectratio]{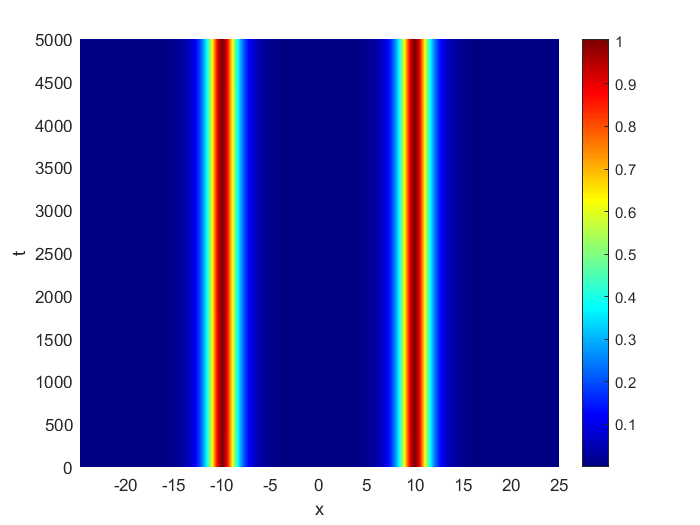}}
\centerline{\text{(a) Periodic Boundary}}
\end{minipage}
\begin{minipage}[b]{0.33\linewidth}
\centering
\centerline{\includegraphics[width=\linewidth,height=\textheight,keepaspectratio]{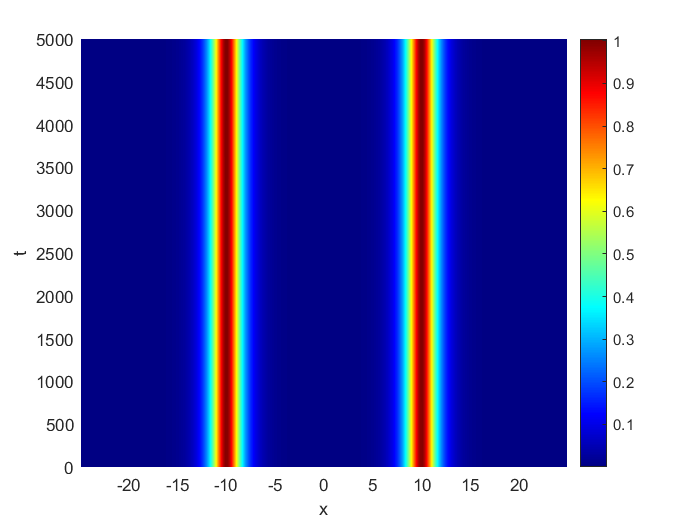}}
\centerline{\text{(b) Neumann Boundary}}
\end{minipage}
\begin{minipage}[b]{0.33\linewidth}
\centering
\centerline{\includegraphics[width=\linewidth,height=\textheight,keepaspectratio]{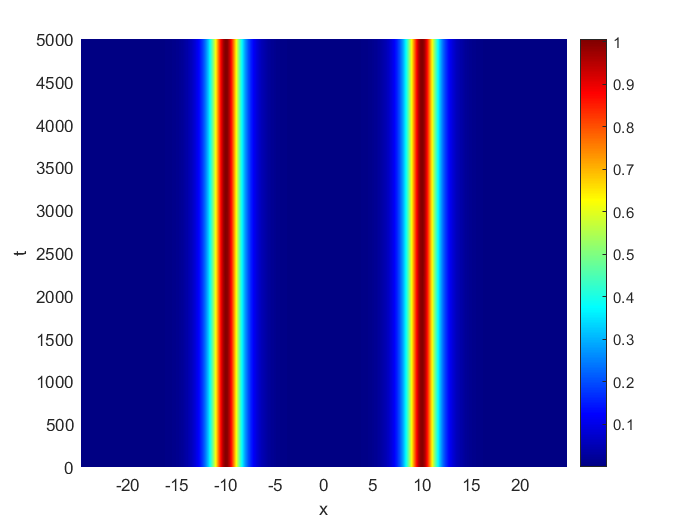}}
\centerline{\text{(c) Dirichlet Boundary}}
\end{minipage}
\caption{\footnotesize{Empirical blow-up analysis of the IFRK4-P13 method}}
\label{fig: ifrkblowup}
\end{figure}

As discussed in \cite{BK, SV}, to evaluate the applicability of the Krogstad-P22 and IFRK4-P13 methods for simulating system \eqref{eq:1} when \( M = 1 \) over long time periods, it is essential to assess whether these methods accurately integrate the spatially independent version of system \eqref{eq:1}. If either method fails to accurately integrate the spatially independent version of system \eqref{eq:1}, it will be prone to nonlinear blowup over a long period. To consider an $x$-independent version of system \eqref{eq:1}, we consider the system \eqref{eq:1} with M=1 and $\alpha_1$=0: 
\begin{flalign}
    &i\frac{\partial \Psi_1}{\partial t} + 2|\Psi_1|^2\Psi_1=0 &
\label{eq:blowup}
\end{flalign}
with initial conditions:
\begin{flalign}
    &\Psi_1(x,0) = \exp{(2i(x+10))}\text{sech}(x+10) + \exp{(-2i(x-10))}\text{sech}(x-10) &
\label{eq:blowupconditions}
\end{flalign}
under periodic, homogeneous Neumann, and homogeneous Dirichlet boundary conditions, with parameters $N=256$ and $k=0.1$ over time $T=5000$. The initial condition models two solitons, and the simulation results utilizing the Krogstad-P22 method under the mentioned boundary conditions are shown in Figure \ref{fig: fig17}. The simulation results utilizing the IFRK4-P13 method under the same boundary conditions are shown in Figure \ref{fig: ifrkblowup}.

From both Figure \ref{fig: fig17} and \ref{fig: ifrkblowup}, we observed stable soliton behavior throughout simulations, with solitons maintaining their forms and velocities. While there is not yet a formal analysis of blow-up that we can consider further, it stands to reason that both Krogstad-P22 and IFRK4-P13 provide good approximations of a system of coupled nonlinear Schrödinger equations.

\subsection{2D CNLSE}
\noindent{\bf Example 4 (Interaction of four waves):}
\\ \indent To check the applicability of the proposed methods for simulating the 2D system \eqref{eq:1} with $M=4$, we modified Example 3 in \cite{Li} to a four-wave interaction by adding two new waves with different initial positions. We also multiply the term $\exp{-i(x^2+y^2)}$ onto the initial conditions of all waves to establish wave movement that induces wave collisions so we may consider their interaction.

Here, we considered the four coupled 2D system \eqref{eq:1} on the domain $\Omega=[-10,10]^2$ with homogeneous Dirichlet boundary conditions and following initial conditions:

\begin{flalign}
 &\Psi_1(x,y,0)=\frac{2}{\sqrt{\pi}}\exp{(-[(x-c)^2+(y-c)^2])}\exp{(-i(x^2+y^2))},&\notag\\
    &\Psi_2(x,y,0)=\frac{2}{\sqrt{\pi}}\exp{(-[(x+c)^2+(y+c)^2])}\exp{(-i(x^2+y^2))},&\notag\\
    &\Psi_3(x,y,0)=\frac{2}{\sqrt{\pi}}\exp{(-[(x-c)^2+(y+c)^2])}\exp{(-i(x^2+y^2))},&\\
    &\Psi_4(x,y,0)=\frac{2}{\sqrt{\pi}}\exp{(-[(x+c)^2+(y-c)^2])}]\exp{(-i(x^2+y^2))}.& \notag
\end{flalign}
To conduct the numerical experiments in Example 4, we considered parameters $\alpha_j=1$, $\sigma_{jj}=1$ for $j=1,2,3,4$, and $\sigma_{jm}=3$ for $j,m=1,2,3,4$ such that $j\neq m$. 

As in previous examples, we performed the temporal convergence of the proposed methods and stored results in Table \ref{tab:Table 7}. In this test, we used \( N = 128 \) and an initial time step \( k = 0.01 \), running the simulation up to \( T = 1 \). As shown in Table \ref{tab:Table 7}, both methods achieved fourth-order accuracy in 2D, as anticipated. 

\begin{figure}[H]
\begin{minipage}[b]{0.49\linewidth} 
\centering
\centerline{\includegraphics[width=\linewidth,height=\textheight,keepaspectratio]{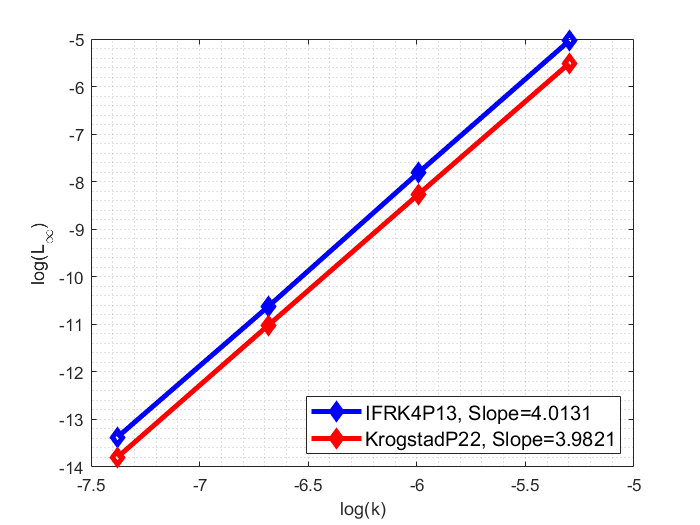}}
\centerline{\text{(a)}}
\end{minipage}
\begin{minipage}[b]{0.49\linewidth}
\centering
\centerline{\includegraphics[width=\linewidth,height=\textheight,keepaspectratio]{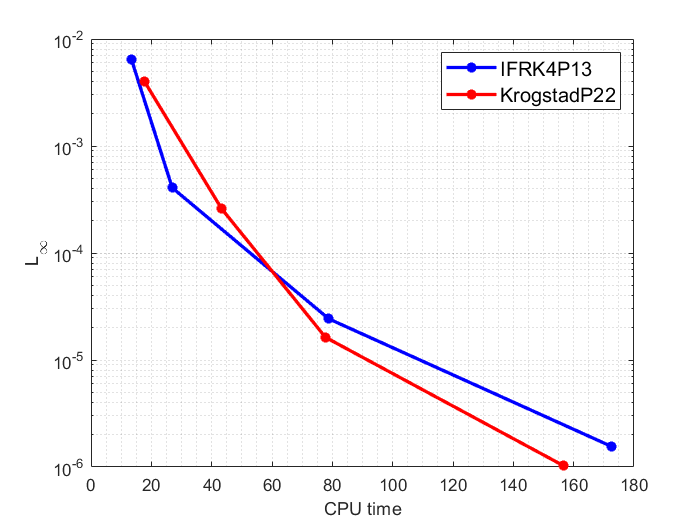}}
\centerline{\text{(b)}}
\end{minipage}
\caption{\footnotesize{Temporal convergence (a) and CPU efficiency (b) for decreasing time steps $k$}}
\label{fig: 2Dconverg}
\end{figure}

\begin{table}[H]
\begin{center} \footnotesize
\resizebox{\linewidth}{!}{
\begin{tabular}{ccccc|ccc}
 \hline
  &\multicolumn{3}{c}{Krogstad-P22}&&\multicolumn{3}{c}{IFRK4-P13}
  \\\cmidrule{2-4}\cmidrule{5-8}
 $k$ & $\norm{\Psi_1}_\infty$ & Order &CPU-time && $\norm{\Psi_1}_\infty$ & Order &CPU-time\\
  \hline\\
  $\frac{1}{200}$& 4.0319$\times 10^{-3}$ & - & 17.671 && 6.4885$\times 10^{-3}$ & - & 13.371 \\\\
  $\frac{1}{400}$& 2.5832$\times 10^{-4}$ & 3.9642 & 43.334 && 4.0670$\times 10^{-4}$ & 3.9958 & 26.942 \\\\
  $\frac{1}{800}$& 1.6311$\times 10^{-5}$ & 3.9853 & 77.706 && 2.4235$\times 10^{-5}$ & 4.0688 & 78.832 \\\\
  $\frac{1}{1600}$& 1.0224$\times 10^{-6}$ & 3.9957 & 156.778 && 1.5615$\times 10^{-6}$ & 3.9561 & 172.438 \\\\
\hline
\end{tabular}}
\end{center} 
\caption{Comparison of accuracy for Krogstad-P22 and IFRK4-P13 methods for Example 4}
\label{tab:Table 7} 
\end{table} 
The temporal convergence data presented in Table \ref{tab:Table 7} is shown on a log-log scale in Figure \ref{fig: 2Dconverg} (a), and the CPU efficiency is shown in Figure \ref{fig: 2Dconverg} (b). From Figure \ref{fig: 2Dconverg} one can see that the Krogstad-P22 method outperforms the IFRK4-P13 method in terms of accuracy and efficiency in 2D as well.

In Figure \ref{fig: 2Dconverg} (b), we compared the computational efficiency of the proposed methods. We observe that both methods approximate with close error and CPU time needed in this 2D example, with Krogstad-P22 only slightly outperforming IFRK4-P13 for smaller time steps.

\begin{figure}[H]
\centering
\centerline{\includegraphics[width=0.49\linewidth,keepaspectratio]{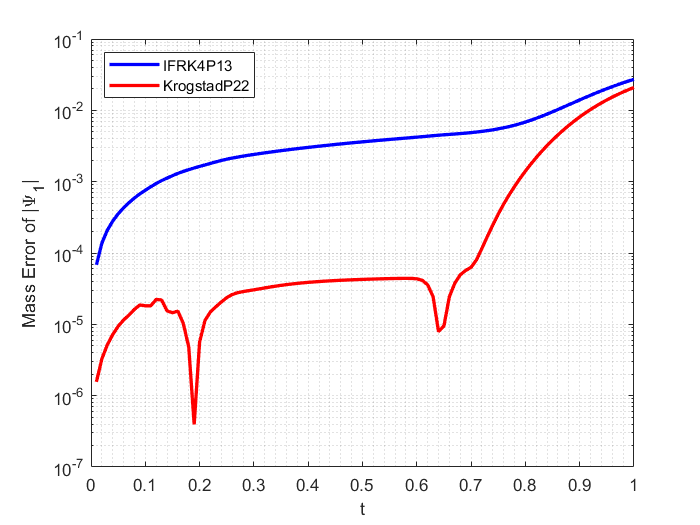}}
\caption{\text{Mass conservation of $\Psi_1$}}
\label{fig: 2Dmass}
\end{figure}

We compared the mass conservation errors of $|\Psi_1|$ for the Krogstad-P22 and IFRK4-P13 methods in Figure \ref{fig: 2Dmass}, with $k=0.01$ and N=128. From Figure \ref{fig: 2Dmass}, we can see that while the methods conserve mass well initially, the mass conservation error grows following the wave interaction. We observed that the Krogstad-P22 method maintains better mass conservation than the IFRK4-P13 method. 

\begin{figure}[H]
\begin{minipage}[b]{0.5\linewidth} 
\centering
\centerline{\includegraphics[width=\linewidth,height=\textheight,keepaspectratio]{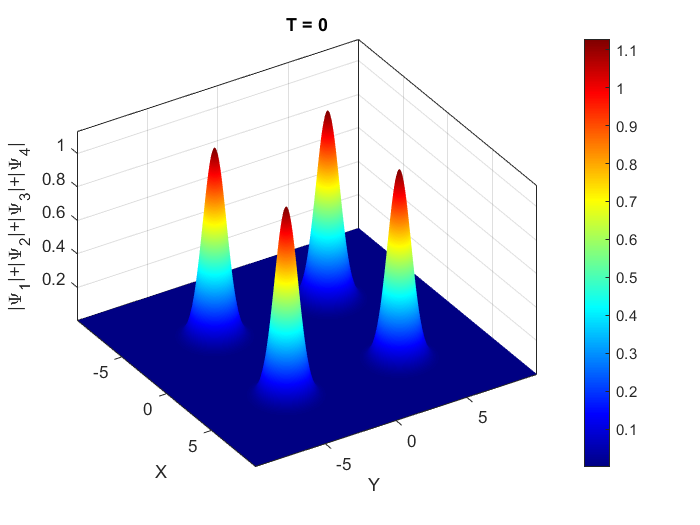}}
\end{minipage}
\begin{minipage}[b]{0.5\linewidth}
\centering
\centerline{\includegraphics[width=\linewidth,height=\textheight,keepaspectratio]{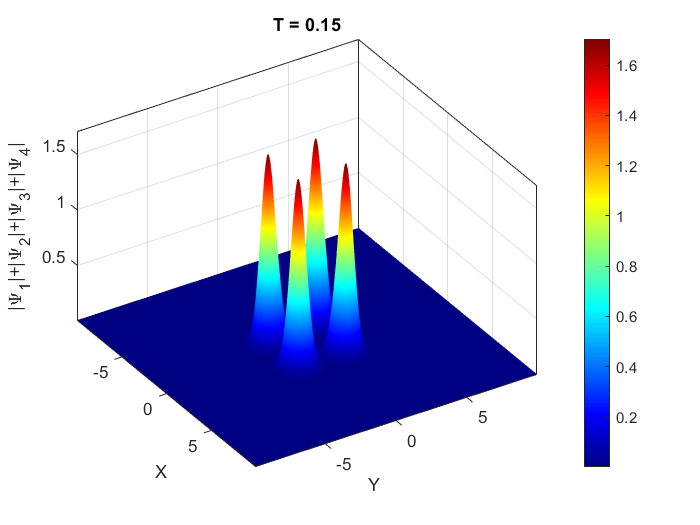}}
\end{minipage}
\begin{minipage}[b]{0.5\linewidth}
\centering
\centerline{\includegraphics[width=\linewidth,height=\textheight,keepaspectratio]{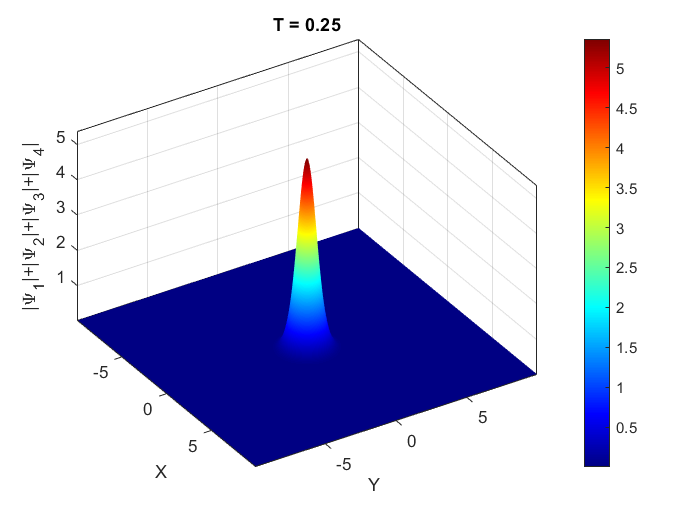}}
\end{minipage}
\begin{minipage}[b]{0.5\linewidth}
\centering
\centerline{\includegraphics[width=\linewidth,height=\textheight,keepaspectratio]{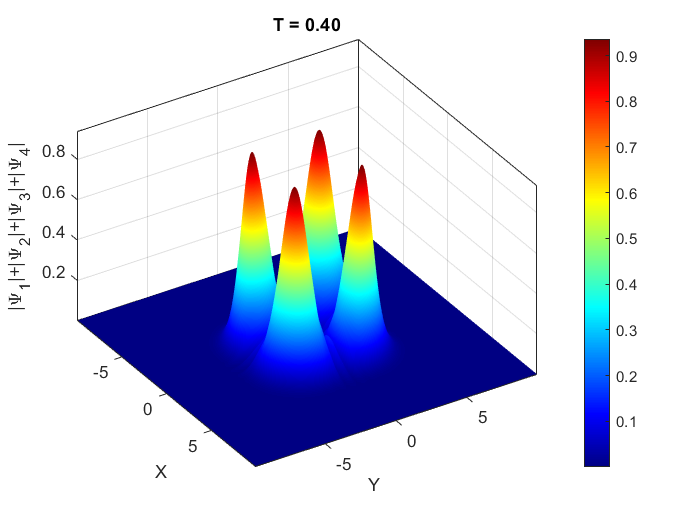}}
\end{minipage}
\caption{\footnotesize{Numerical simulation of 2D wave interaction at various times $T=0,0.15,0.25,0.40$.}}
\label{fig: 2Dint}
\end{figure}

We illustrated the simulation of the nonlinear wave interaction for this example. The wave interaction profiles are shown in Figure \ref{fig: 2Dint}. From the simulated solutions at $T = 0, 0.15, 0.25, 0.40$, we observe the initial conditions, and the waves before, during, and after the wave interaction respectively. While the waves separate after interaction without the branching or fusion of any waves, $|\Psi(x,y,t)|$ has a lower value for each wave than initially observed, suggesting this interaction is inelastic. We will next consider a 3D extension of this scenario to investigate this simulation further.

\subsection{3D CNLSE}
\noindent{\bf Example 5 (Interaction of four waves):}
\\ For our final instance, we consider the 3D extension of the scenario in Example 4 to further analyze the use of both methods in higher dimensions. We consider system \eqref{eq:1} with $M=4$ on the domain $\Omega=[-10,10]^3$ to model the interaction of four waves in a nonlinear medium. The initial conditions that we consider are as follows, extending those from Example 4 to be 3D:
\begin{flalign}
    &\Psi_1(x,y,z,0)=\frac{2}{\sqrt{\pi}}\exp{(-[(x-c)^2+(y-c)^2+(z-c)^2])}\exp{(-i(x^2+y^2+z^2))}, \notag&\\\vspace{0.3cm}  
    &\Psi_2(x,y,z,0)=\frac{2}{\sqrt{\pi}}\exp{(-[(x+c)^2+(y+c)^2+(z+c)^2])}\exp{(-i(x^2+y^2+z^2))}, \notag&\\\vspace{0.3cm}
    &\Psi_3(x,y,z,0)=\frac{2}{\sqrt{\pi}}\exp{(-[(x-c)^2+(y+c)^2+(z-c)^2])}\exp{(-i(x^2+y^2+z^2))}, \notag&\\\vspace{0.3cm}
    &\Psi_4(x,y,z,0)=\frac{2}{\sqrt{\pi}}\exp{(-[(x+c)^2+(y-c)^2+(z+c)^2])}\exp{(-i(x^2+y^2+z^2))}\notag.&
\end{flalign}
with homogeneous Dirichlet boundary conditions and parameters $\alpha_j=1$, $\sigma_{jj}=1$ for $j=1,2,3,4$, and $\sigma_{jm}=3$ for $j,m=1,2,3,4$ such that $j\neq m$.

\begin{figure}[H]
\begin{minipage}[b]{0.49\linewidth} 
\centering
\centerline{\includegraphics[width=\linewidth,height=\textheight,keepaspectratio]{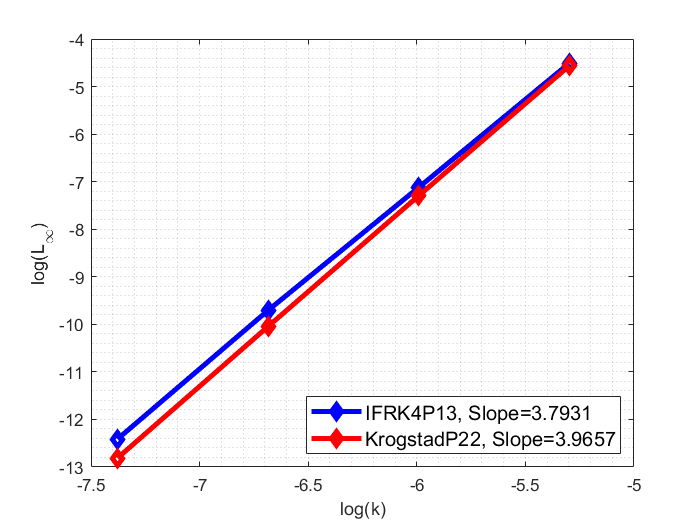}}
\centerline{\text{(a)}}
\end{minipage}
\begin{minipage}[b]{0.49\linewidth}
\centering
\centerline{\includegraphics[width=\linewidth,height=\textheight,keepaspectratio]{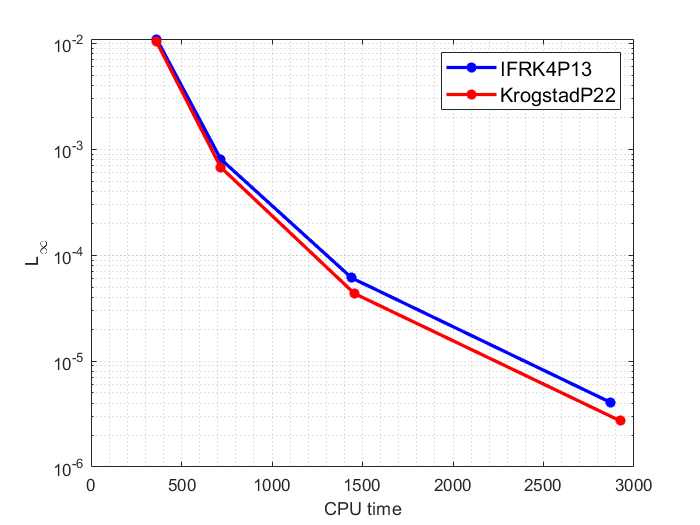}}
\centerline{\text{(b)}}
\end{minipage}
\caption{\footnotesize{Pictorial representation of Table \ref{tab:Table 8} for temporal convergence (a) and computational efficiency (b) }}
\label{fig: 3Dconverg}
\end{figure}

To compare the performance of both methods for simulating this 3D problem, we analyzed the temporal convergence of the methods for halving time step $k=0.01$ keeping $N=64$ fixed up to $T=1$, with the results shown in Table \ref{tab:Table 8}. From Table \ref{tab:Table 8} and Figure \ref{fig: 3Dconverg} (a), it is evident that the Krogstad-P22 method achieved fourth-order convergence in the temporal direction, although the IFRK4-P13 method struggles to do so. From Figure \ref{fig: 3Dconverg} (b), we can notice that the Krogstad-P22 method is slightly more computationally efficient than the IFRK4-P13 method in the 3D case. These results suggest that the Krogstad-P22 method may be better suited for simulating nonlinear PDEs in 3D.
\begin{table}[H]
\begin{center} \footnotesize
\resizebox{\linewidth}{!}{
\begin{tabular}{ccccc|ccc}
 \hline
  &\multicolumn{3}{c}{Krogstad-P22}&&\multicolumn{3}{c}{IFRK4-P13}
  \\\cmidrule{2-4}\cmidrule{5-8}
 $k$ & $\norm{\Psi_1}_\infty$ & Order &CPU-time && $\norm{\Psi_1}_\infty$ & Order &CPU-time\\
  \hline\\
  $\frac{1}{200}$& 1.0418$\times 10^{-2}$ & - & 360.935 && 1.1002$\times 10^{-2}$ & - & 362.584 \\\\
  $\frac{1}{400}$& 6.7796$\times 10^{-4}$ & 3.9418 & 716.325 && 8.0750$\times 10^{-4}$ & 3.7681 & 716.668 \\\\
  $\frac{1}{800}$& 4.3438$\times 10^{-5}$ & 3.9642 & 1455.434 && 6.0925$\times 10^{-5}$ & 3.7284 & 1440.680 \\\\
  $\frac{1}{1600}$& 2.7306$\times 10^{-6}$ & 3.9917 & 2923.362 && 4.0690$\times 10^{-6}$ & 3.9043 & 2870.947 \\\\
\hline
\end{tabular}}
\end{center} 
\caption{Comparison of temporal convergence for the Krogstad-P22 and IFRK4-P13 methods as used in Example 4}
\label{tab:Table 8} 
\end{table} 

\begin{figure}[H] 
\centering
\centerline{\includegraphics[width=0.49\linewidth,keepaspectratio]{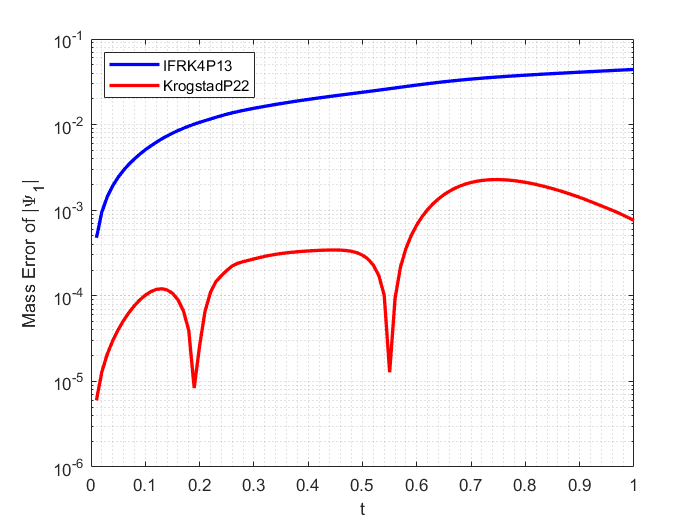}}
\caption{\text{Mass conservation of $\Psi_1$}}
\label{fig: 3Dmass}
\end{figure}

In Figure \ref{fig: 3Dmass}, we analyze the mass conservation error of $|\Psi_1|$. The mass approximation errors seen in Figures \ref{fig: 2Dmass} and \ref{fig: 3Dmass} suggested that the methods compared may struggle with mass conservation over long durations in higher-dimensional simulations, potentially due to inelastic interactions. However, the Krogstad-P22 method preserved wave mass accurately up to three decimal places over short periods, allowing for reliable 3D wave interaction simulations for up to $T=1$ before significant mass approximation errors occur.

With this in mind, we simulated this example up to $T=1$. We use the Krogstad-P22 method on a $N \times N \times N$ grid, resulting in a grid of wave amplitude value. To visualize this result properly, we use MATLAB to find an isovlaue, which we then use to create an isosurface from our grid. Below in Figure \ref{fig: 3Dint} we can see the simulation before and during wave interaction.

\begin{figure}[H]
\begin{minipage}[b]{0.5\linewidth} 
\centering
\centerline{\includegraphics[width=\linewidth,height=\textheight,keepaspectratio]{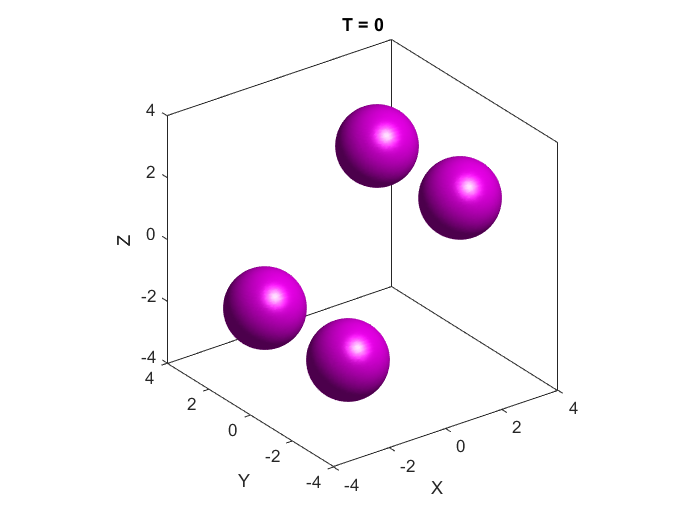}}
\end{minipage}
\begin{minipage}[b]{0.5\linewidth}
\centering
\centerline{\includegraphics[width=\linewidth,height=\textheight,keepaspectratio]{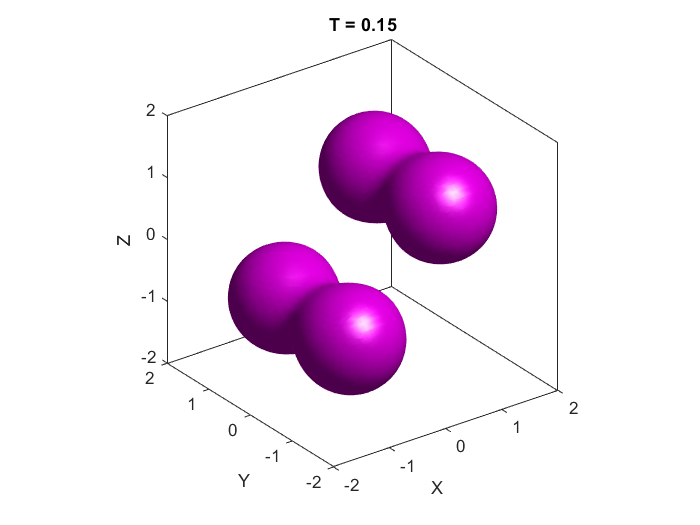}}
\end{minipage}
\begin{minipage}[b]{0.5\linewidth}
\centering
\centerline{\includegraphics[width=\linewidth,height=\textheight,keepaspectratio]{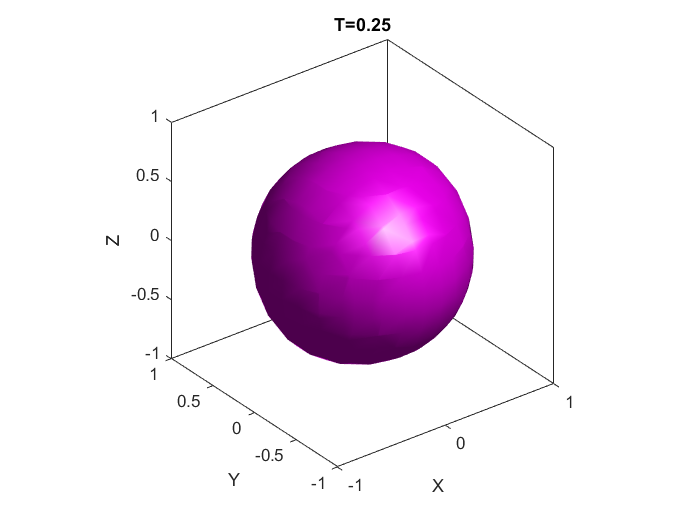}}
\end{minipage}
\begin{minipage}[b]{0.5\linewidth}
\centering
\centerline{\includegraphics[width=\linewidth,height=\textheight,keepaspectratio]{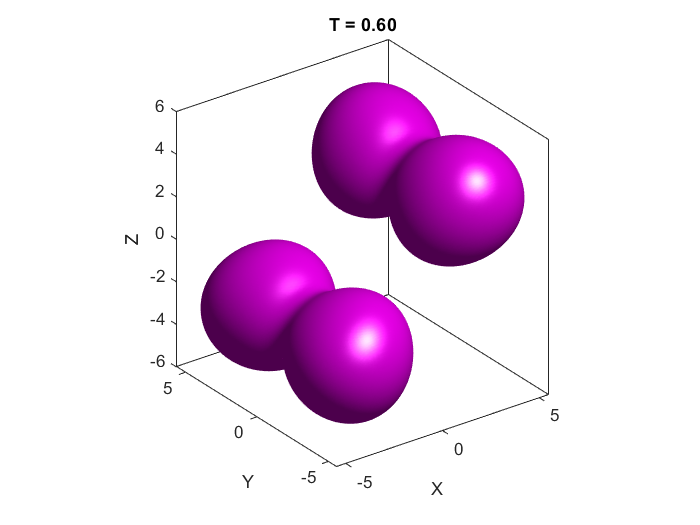}}
\end{minipage}
\caption{\footnotesize{Spatio-temporal wave propagation isosurface of $|\Psi_1|+|\Psi_2|+|\Psi_3|+|\Psi_4|$ at various times $T=0,0.15,0.25$, and $0.60$.}}
\label{fig: 3Dint}
\end{figure}

\section{Conclusions}
\label{sec:5}
This paper examined the effectiveness of the Krogstad-P22 and IFRK4-P13 methods, combined with the Fourier spectral method, for simulating multi-dimensional M-coupled nonlinear Schrödinger equations. Both methods achieved fourth-order temporal accuracy and spectral spatial convergence in applicable cases, with the exception of the IFRK4-P13 method in Example 4. The Krogstad-P22 method demonstrated superior computational efficiency and accuracy, achieving lower errors in shorter or equivalent CPU times compared to IFRK4-P13, across one-, two-, and three-dimensional simulations. 

Under various boundary conditions, including periodic, homogeneous Neumann, and homogeneous Dirichlet, the Krogstad-P22 method maintained better accuracy and stability. It also preserved mass and energy conservation properties more effectively, highlighting its potential for long-term simulations. Future research may extend these methods to more complex systems and explore applications in other areas requiring precise multi-dimensional nonlinear dynamics simulations.

\section{Acknowledgments}
This research was funded by the Turner Endowment for Engaged Learning in STEM (TEELS) grant through Utah Valley University.

\bibliographystyle{unsrt}  



\end{document}